\documentclass{article}
\usepackage[utf8]{inputenc}
\usepackage[english]{babel}
\usepackage{fullpage,amsmath,amssymb,amsthm,amsfonts,enumerate,textcomp,eurosym,titling,epsfig,epstopdf,tikz,esint,bigints,mathrsfs,todonotes,mathtools,verbatim,bm,nccmath,csquotes}
\usepackage{tikz-cd}
\usepackage{dirtytalk}
\usepackage{titlesec}
\usepackage[style=numeric]{biblatex}
\input xy
\xyoption{all}

\titleformat{\paragraph}
{\normalfont\normalsize\bfseries}{\theparagraph}{1em}{}
\titlespacing*{\paragraph}
{0pt}{3.25ex plus 1ex minus .2ex}{1.5ex plus .2ex}

\DeclareMathOperator{\Eu}{Eu}

\numberwithin{equation}{section}

\title{Relations in Twisted Quantum $K$-Rings}
\author{Irit Huq-Kuruvilla}
\begin{document}

\maketitle

\theoremstyle{plain}
\newtheorem{thm}{Theorem}[section]
\newtheorem{lem}[thm]{Lemma}
\newtheorem{prop}[thm]{Proposition}
\newtheorem{cor}[thm]{Corollary}
\newtheorem{exercise}[thm]{Exercise}

\newtheorem*{thm*}{Theorem}

\theoremstyle{definition}
\newtheorem{mydef}[thm]{Definition}
\newtheorem{conj}[thm]{Conjecture}
\newtheorem*{conj*}{Conjecture}
\newtheorem{exmp}[thm]{Example}
\newtheorem*{claim}{Claim}
\newtheorem*{notation}{Notation}
\newtheorem{prob}{Problem}
\newtheorem{ex}{Exercise}

\theoremstyle{remark}
\newtheorem*{rem}{Remark}
\newtheorem*{note}{Note}

\newcommand{\Z}{\mathbb{Z}}
\newcommand{\Q}{\mathbb{Q}}
\newcommand{\R}{\mathbb{R}}
\newcommand{\C}{\mathbb{C}}
\newcommand{\bb}{\mathbb}
\newcommand{\cali}{\mathcal}
\newcommand{\oo}{\omega}
\newcommand{\p}{\partial}
\newcommand{\Wk}{W^{k,p}(U)}
\newcommand{\W}{W^{1,p}(U)}
\newcommand{\g}{\mathfrak{g}}
\newcommand{\fk}{\mathfrak}
\newcommand{\ft}{\text{ft}}
\newcommand{\ev}{\operatorname{ev}}
\newcommand{\h}{\mathfrak{h}}
\newcommand{\mbar}{\overline{M}}
\newcommand{\cO}{\mathcal{O}}
\newcommand{\J}{\mathcal{J}}
\newcommand{\Kloop}{\mathcal{K}}
\newcommand{\llangle}{\left\langle\!\left\langle}
\newcommand{\rrangle}{\right\rangle\!\right\rangle}
\newcommand{\sslash}{\mathbin{/\mkern-6mu/}}
\newcommand{\gl}{\operatorname{gl}}
\newcommand{\st}{\operatorname{st}}

\begin{abstract}
We introduce twisted quantum $K$-rings, defined via inserting a twisting class into the virtual structure sheaf. We develop a toolkit for computing relations by adapting some results about ordinary quantum K rings to our setting, and discuss some applications, including Ruan-Zhang's quantum $K$-theory with level structure and complete intersections inside $\mathbb{P}^n$,  proving some conjectures obtained physically via 3D gauged linear sigma models. We also twisted version of Lee-Pandharipande's reconstruction theorem. 

In addition, we formulate a ring-theoretic  abelian/non-abelian correspondence conjecture, relating the quantum K-ring of a GIT quotient $V//G$ to a certain twist of the quantum K-ring of $V//T$, the quotient by the maximal torus. We prove this conjecture for partial flags, and use this to give another proof of the Whitney presentation for the quantum $K$-theory ring of the Grassmanian.

 We generalize all the above work to an $S_n$-equivariant setting.

\end{abstract}

\section*{Introduction}

The quantum $K$-ring (introduced by Givental and Lee) of a smooth projective variety $X$ is a deformation of its topological $K$-ring involving $K$-theoretic Gromov-Witten invariants of $X$. These invariants are defined as certain holomorphic Euler characteristics of sheaves on Kontsevich's moduli space of stable maps $\overline{M}_{0,n,d}(X)$. 

It is a $K$-theoretic counterpart to the quantum cohomology ring introduced by Witten and Kontsevich. Many recent developments in physics (including those of \cite{jock}, \cite{ueda}) have resulted in connections to quantum $K$-theory, and predictions for relations in quantum $K$-rings. However, quantum $K$-rings are in general much harder to compute than their cohomological counterparts, so many physical predictions go untested.

In addition, much work in the subject of quantum $K$-theory has been about the quasimap $K$-rings of Nakajima quiver varieties, and their relationships with quantum groups, defined with respect to a different moduli space. A natural question to ask is how these rings are related to their \say{genuine} stable-map counterparts. 

The goal of this work is to develop theoretical groundwork for answering some of these questions by developing the theory of twisted quantum $K$-rings, based on certain modifications of the virtual structure sheaf, and to give some immediate applications of the resulting theorems to complete intersections, cotangent bundles, and GIT quotients. (Note the word twisted is used as in \cite{cg}, these are not orbifold Gromov-Witten invariants). 

It was observed by Givental that one can calculate Gromov-Witten invariants of the quintic (denoted $Q_5$) by doing intersection theory on $\overline{M}_{0,n,d}(\mathbb{P}^4)$, with the virtual fundamental class modified by the Euler class of the normal bundle $\overline{M}_{0,n,d}(Q_5)$, which can be expressed as a pushforward from $\overline{M}_{0,n+1,d}$ of $ev_{n+1}^*\mathcal{O}(5)$.   This introduced the theory of twisted Gromov-Witten invariants, fully developed in \cite{cg}. In general, twisted invariants can be computed from their untwisted counterparts. 

The $K$-theoretic counterpart of twisted invariants were introduced in \cite{tonitaktwist}, and generalizations were developed in \cite{PermXI} and \cite{QuantCob2}. In this paper, we will consider the most common type of twisting (type I), and develop a theory of quantum $K$-rings based on these invariants.

Fixing $g,n,d$, we have the diagram
\[
\xymatrix{\overline{M}_{g,n+1,d}(X)\ar[d]^{\ft}\ar[r]^(.65){\ev_{n+1}} & X\\ \overline{M}_{g,n,d}(X) &}.
\]
Here $\overline{M}_{g,n+1,d}(X)$ is the universal family of
$\overline{M}_{g,n,d}(X)$, $\ft$ forgets the last marked point, and
$\ev_{n+1}$ evaluates at that point.

Essentially, a type-I twisting is a constraint pulled back by $\ev_{n+1}$
and then pushed forward by $\ft$, which imposes a condition on the entire
map rather than only at one marked point.

More precisely, given any $V\in K(X)$ and such a multiplicative
$K$-theoretic characteristic class $C$, twisted $K$-theoretic
Gromov--Witten invariants with twisting $C(V)$ are holomorphic Euler
characteristics on $\overline{M}_{g,n,d}(X)$ with respect to
$\mathcal O^{tw}:=\mathcal O^{vir}\otimes C(\ft_*\ev_{n+1}^*V)$.

As it turns out, type-I twisted invariants satisfy very similar properties to their untwisted counterparts, and as such, many results on quantum $K$-rings pass unchanged to the twisted case, in particular, we show the following:

\begin{thm*}
    Type I twisted invariants satisfy analogues of the Kontsevich-Manin axioms, and they are the structure constants to a ring, the (big or small) \emph{twisted quantum $K$ ring} of $X$.
\end{thm*}

Furthermore we introduce the twisted $J$ function $\mathcal{J}^{tw}$, a certain kind of generating function for twisted invariants, and we prove the following theorem (which was proven by Iritani-Milanov-Tonita in \cite{imt} in the untwisted case):

\begin{thm*}
    If $K(X)$ is generated by line bundles whose first Chern classes are effective curve classes, then
    a polynomial difference operator in $q$ annihilating $\mathcal{J}^{tw}$ yields a relation in the twisted quantum $K$-ring of $X$.
\end{thm*}

Using the above theorem in general is difficult, as the relations are obtained in terms of operators $A_{i}$ which arise as solutions to a Lax-type equation. We give some conditions on which these operators can be simply interpreted in terms of multiplication by classes in $K(X)$, generalizing some results of Anderson-Chen-Tseng-Iritani in \cite{act}.  We work in the generality of $S_n$-equivariant quantum $K$-theory, as developed by Givental.

As examples, we use the results to calculate relations in the (small) quantum $K$-rings for complete intersections in $\mathbb{P}^N$, as well as Ruan-Zhang's quantum $K$-theory with level structure, and we show that the quasimap quantum $K$-ring for $T^*\mathbb{P}^N$ introduced in \cite{Bax} is equal to a specific modification of its stable map counterpart.

We also introduce an abelian/non-abelian correspondence conjecture for small quantum $K$-rings. 
Our setting will be some reductive group $G$ acting on a vector space $V$. The classical abelian/non-abelian correspondence proven by Martin constructs a ring map $\phi:H^*(V//T)^W\to H^*(V//G)$, and identifies the Poincare pairing on the target with the Poincare pairing on the source modified by the Euler class of the vector bundle $R^M:=\sum_{\alpha \in \text{roots($G$)}} L_\alpha$. The analogous statement in $K$-theory was proven by Harada-Landweber \cite{harada}. We propose the following quantum version of that conjecture.
\begin{conj*}
    The map $x\mapsto \lim_{\lambda\to 1} \phi(x)$, extended to Novikov variables, is a surjective ring homomorphism from $QK^{tw}(V\sslash T)$ to $QK(V\sslash G)$, where the twisting is determined by the class $\Eu_\lambda(R_M)$, where $Eu_\lambda$ is the $\mathbb{C}^*$-equivariant Euler class, where $ \mathbb{C}^*$ scales the fibers of $R$. 
\end{conj*}

We note that this is not a $K$-theoretic version of the correspondence for quantum cohomology rings conjectured by Ciocan-Fontanine--Kim--Sabbah in \cite{cfks}, that conjecture relates untwisted rings. There is a likely a $K$-theoretic counterpart to that conjecture, however the version we work with matches up well with the kind of formulas that appear in physics for relations in quantum $K$-rings, as we will see for the Grassmanian.
We show this conjecture is a consequence of other abelian/nonabelian-type conjectures for $J$-functions, and, using the theorem of Givental--Yan \cite{Xiaohan} on the symmetrized $J$-function of the Grassmannian, we prove the conjecture for Grassmannians. Using this, we reprove the presentation for $QK_T(\operatorname{Gr}(k,n))$ from \cite{gu2022quantum}. Our arguments should extend to all partial flag varieties, but we leave that case to future work.

\section{Twisted $K$-Theoretic Gromov-Witten Invariants}
\label{defs}
\subsection{Basic Setup}
Let $X$ be a smooth projective variety, possibly equipped with an action of an
algebraic torus $T$.  Throughout the paper, $K(X)$ denotes complex
topological $K$-theory and $K_T(X)$ denotes $T$-equivariant complex
topological $K$-theory. 

$K$-theoretic Gromov-Witten invariants are ($T$-equivariant) holomorphic Euler characteristics of sheaves on the space $\overline{M}_{g,m}(X,d)$, Kontsevich's moduli space of degree $d$ stable maps from $n$-pointed genus $g$ curves to $X$, henceforth denoted $X_{g,m,d}$, taken a gainst the virtual structure sheaf $\mathcal{O}^{vir}$.

We work over a $\lambda$-algebra $\Lambda$ with Adams operations $\Psi^r$,
$r\geq1$, containing the Novikov variables.  Let
$K_T(X)\otimes\Lambda$ denote the completed tensor product.  We fix an ideal
$\Lambda_+\subset K_T(X)\otimes\Lambda$ containing every positive-degree
Novikov monomial and assume that $K_T(X)\otimes\Lambda$ is complete with
respect to $\Lambda_+$.  The Adams operations extend to
$K_T(X)\otimes\Lambda$, and we impose $\Psi^r(Q^d)=Q^{rd}$ on Novikov
monomials. In the presence of nontrivial torus action, we let the torus parameters in $K_T(pt)$ be denote $\Lambda_i$, and we enforce $\Lambda_i\in \Lambda$, and require $\Lambda$ is completed with respect to each $1-\Lambda_i$. The default example of such a setup is $\Lambda=\hat{K}_T(pt)[Q]$, where the hat denotes completion at $Q$ and $1-\Lambda_i$. So the ideal $\Lambda_+\in K_T(X)\otimes \Lambda$ is the ideal generated by $1-\Lambda_i$ and the Novikov variables.

\subsection{Type-I Twistings}

A multiplicative characteristic class is a monoid homomorphism
$C:(\operatorname{Vect},\oplus)\to(K_T(-)\otimes\Lambda)^\times$.
An invertible such class extends to virtual bundles by
$C(A-B)=C(A)/C(B)$.  We will use classes of the form
$C(E)=\exp\bigl(\sum_{r\neq0}s_r\Psi^r(E)\bigr)$, with $s_r\in\Lambda$, in
a topology in which the expression converges.

Let $\alpha$ be an element of $K_T(X_{g,m+1,d})$.  Define $E_{g,m,d}(\alpha):=R\ft_*\alpha.$ A \emph{twisting} replaces the virtual structure sheaf
$\cO^{\mathrm{vir}}_{g,m,d}$ with $\cO^{\mathrm{vir},\mathrm{tw}}_{g,m,d}
:=\cO^{\mathrm{vir}}_{g,m,d}\otimes C(E_{g,m,d}(\alpha)).$

Twistings are classified into types depending on the choice of $K$-theory element $\alpha$. 
\begin{itemize}
\item Type I uses $\alpha=\ev_{m+1}^*V$, for
$V\in K_T(X)\otimes\Lambda$, as introduced in
\cite[\S4]{tonitaktwist}. This can be conceptualized as choosing a universal constraint on the image of the stable map in the target. 
\item Type II uses $\alpha=\ev_{m+1}^*f(L_{m+1})$, where $f$ is a Laurent
polynomial with coefficients in $K_T(X)$ and $f(1)=0$, as introduced in
\cite[Theorem~2]{PermXI}.  A general Laurent polynomial $f$ can be separated into type-I and type-II
parts.
\item Type III uses $\alpha=i_*f(L_+,L_-)$ on the nodal locus of the
universal curve, where $f$ is symmetric in the cotangent lines to the two
branches, as introduced in \cite[Theorem~3.3]{QuantCob2}.
\end{itemize}
Successive twistings are allowed, and quantum Adams--Riemann--Roch expresses
the resulting permutation-equivariant theory in terms of the untwisted one
\cite[Theorems~1--2]{PermXI}; the nodal transformation is
\cite[Theorem~3.3]{QuantCob2}.

We restrict ourselves to Type I theories, as the others in general do not define rings.

\subsection{$S_n$-equivariant invariants}

We use Givental's definition of mixed permutation-equivariant invariants developed in \cite{PermVII}. The group $S_n$ acts on $X_{g,k+n,d}$ by fixing the first $k$ markings and permuting the last $n$ markings. We call the $n$ marked points \emph{permutable marked points}, and the $k$ unpermutable markings \emph{named marked points}. 

\begin{mydef}
\label{mixed-equivariant-correlator-definition}
Let $a_1(q),\dots,a_k(q)$ and $t(q)$ be Laurent polynomials in $q$ with
coefficients in $K_T(X)\otimes\Lambda$.  Write
$a_i(q)=\sum_m a_{i,m}q^m$ and $t(q)=\sum_m t_mq^m$.  Set
$\ev_i^*a_i(L_i):=\sum_mL_i^m\ev_i^*(a_{i,m})$ and
$\ev_{k+j}^*t(L_{k+j}):=\sum_mL_{k+j}^m\ev_{k+j}^*(t_m)$.
Equip the type-I twisted virtual structure sheaf with the natural
$S_n$-action induced by the permutation of the last $n$ markings.
Define the $S_n$-equivariant class
\begin{align*}
\mathscr E^{\mathrm{tw}}_{g,k,n,d}(a;t):={}&
\cO^{\mathrm{vir},\mathrm{tw}}_{g,n+k,d}
\otimes
\bigotimes_{i=1}^k\ev_i^*a_i(L_i)
\otimes
\bigotimes_{j=1}^n\ev_{k+j}^*t(L_{k+j}).
\end{align*}
Here $S_n$ acts simultaneously on the last $n$ markings and on the last
$n$ tensor factors.  Let $\overline{\mathscr E}^{\mathrm{tw}}_{g,k,n,d}(a;t)$
denote the descent of this class to $X_{g,n+k,d}/S_n$.  Define the invariant
to be
\begin{align*}
&\chi_T\left(
 X_{g,n+k,d}/S_n;
 \overline{\mathscr E}^{\mathrm{tw}}_{g,k,n,d}(a;t)
 \right)\notag\\
&=
\sum_{j\geq0}(-1)^j
\left[
 H^j\left(
 X_{g,n+k,d};
 \mathscr E^{\mathrm{tw}}_{g,k,n,d}(a;t)
 \right)
\right]^{S_n}
\in K_T(\mathrm{pt})\otimes\Lambda.
\end{align*}

\end{mydef}

\begin{mydef}
\label{virtual-Sn-module-definition}
The virtual $S_n$-module associated to the insertions in
Definition~\ref{mixed-equivariant-correlator-definition} is
\[
\mathbb V^{\mathrm{tw}}_{g,k,n,d}(a;t)
:=
\sum_{j\geq0}(-1)^j
\left[
 H^j\left(
 X_{g,n+k,d};
 \mathscr E^{\mathrm{tw}}_{g,k,n,d}(a;t)
 \right)
\right]
\in K_T(\mathrm{pt})\otimes\Lambda\otimes R(S_n).
\]
\end{mydef}

With this notation, the invariant in
Definition~\ref{mixed-equivariant-correlator-definition} is
$[\mathbb V^{\mathrm{tw}}_{g,k,n,d}(a;t)]^{S_n}$.

In fact, enlarging the base algebra and coupling the $n$ permutable markings
to representations of the general linear group recovers the entire
$S_n$-module, not merely its invariant part.  This justifies the term
$S_n$-equivariant rather than $S_n$-invariant.

\begin{mydef}
\label{permutation-coefficient-algebra-definition}
Let $\mathrm{Sym}_{\mathbb Q}:=\mathbb Q[p_1,p_2,\dots]$, graded by
$\deg p_r=r$, and let
$\widehat{\mathrm{Sym}}_{\mathbb Q}:=\prod_{\nu\geq0}
\mathrm{Sym}_{\mathbb Q}^{(\nu)}$.
The coefficient algebra used for the stable form of the permutation-equivariant
theory is
 $\Lambda^{\mathrm{perm}}
 :=
 \Lambda\widehat\otimes_{\mathbb Q}
 \widehat{\mathrm{Sym}}_{\mathbb Q}.$
 
Its Adams operations extend those on $\Lambda$ and satisfy $\Psi^r(p_m)=p_{rm}$, in particular $\Psi^r(p_1)=p_r.$

\end{mydef}

\begin{mydef}
\label{rho-N-definition}
Fix $N\geq1$, and let $U_N=\mathbb C^N$ be the standard representation of
$GL_N$.  Define the graded $\mathbb Q$-algebra homomorphism
$\rho_N:\mathrm{Sym}_{\mathbb Q}\to R(GL_N)_{\mathbb Q}$ by
$\rho_N(p_r):=\Psi^r([U_N])$,
where $\Psi^r$ is the $r$th Adams operation in the representation ring.  We
also write $\rho_N$ for the homomorphism acting as the identity on
$K_T(\mathrm{pt})\otimes\Lambda$ and by the preceding formula on the
symmetric-function factor; it acts degree by degree on the completion.  If
$A=\operatorname{diag}(x_1,\dots,x_N)\in GL_N$, then
$\operatorname{Tr}(A\mid\rho_N(p_r))=x_1^r+\cdots+x_N^r$.
Thus $\rho_N$ replaces the formal power sums by the characters of concrete
$GL_N$-representations.
\end{mydef}

\begin{mydef}
\label{expanded-correlator-definition}
Fix $N\geq1$ and include one copy of $U_N$ at each permutable marking before
descending to $X_{g,n+k,d}/S_n$.  Throughout the rest of the paper, every
mixed correlator is understood in the following $T\times GL_N$-equivariant
form, and the fixed parameter $N$ is omitted from its notation:
\begin{align*}
&\left\langle
 a_1,\dots,a_k;t,\dots,t
 \right\rangle^{S_n,\mathrm{tw}}_{g,n+k,d}
 \notag\\
&:=
\chi_{T\times GL_N}\left(
 X_{g,n+k,d}/S_n;
 \overline{
 \mathscr E^{\mathrm{tw}}_{g,k,n,d}(a;t)
 \otimes U_N^{\otimes n}}
 \right)\notag\\
&=
\left[
 \mathbb V^{\mathrm{tw}}_{g,k,n,d}(a;t)
 \otimes U_N^{\otimes n}
\right]^{S_n}.
\end{align*}
Here $S_n$ acts diagonally: it permutes the last $n$ markings, the
corresponding factors of the integrand, and the tensor factors of
$U_N^{\otimes n}$.
After $N$ is fixed, extension of coefficient modules and tensor products by
$R(GL_N)_{\mathbb Q}$ is understood unless an individual
$S_n$-isotypical coefficient is being extracted.
\end{mydef}

\begin{thm}
\label{expanded-invariant-recovery}
Let $S^\lambda$ be the irreducible $S_n$-module indexed by
$\lambda\vdash n$.  Write
$\mathbb V=\mathbb V^{\mathrm{tw}}_{g,k,n,d}(a;t)$ and put
$m_\lambda:=[\operatorname{Hom}_{S_n}(S^\lambda,\mathbb V)]$.
For every $N\geq1$, the correlator of
Definition~\ref{expanded-correlator-definition} has the decomposition
\begin{align}
&\left\langle
 a_1,\dots,a_k;t,\dots,t
 \right\rangle^{S_n,\mathrm{tw}}_{g,n+k,d}
 \notag\\
&=
\sum_{\substack{\lambda\vdash n\\\ell(\lambda)\leq N}}
 m_\lambda\,[\mathbb S_\lambda(U_N)]
\quad\text{in }K_T(\mathrm{pt})\otimes\Lambda\otimes R(GL_N)_{\mathbb Q}.
\label{mixed-Schur-Weyl-extraction}
\end{align}
Moreover, every coefficient $m_\lambda$ is itself the Gromov--Witten
invariant
\begin{align*}
m_\lambda
&=
\chi_T\left(
 X_{g,n+k,d}/S_n;
 \overline{
 \mathscr E^{\mathrm{tw}}_{g,k,n,d}(a;t)
 \otimes(S^\lambda)^\vee}
 \right).
\end{align*}
If $N\geq n$, then every partition of $n$ occurs in
\eqref{mixed-Schur-Weyl-extraction}; hence the coefficients of the irreducible
polynomial $GL_N$-representations $\mathbb S_\lambda(U_N)$ recover every
irreducible multiplicity of the virtual $S_n$-module
$\mathbb V^{\mathrm{tw}}_{g,k,n,d}(a;t)$.  In particular, the invariant of
Definition~\ref{mixed-equivariant-correlator-definition} is the coefficient
$m_{(n)}$, because $S^{(n)}$ is the trivial $S_n$-module.
\end{thm}

\begin{proof}
Schur--Weyl duality gives an equality of $GL_N\times S_n$-modules
\[
 U_N^{\otimes n}
 \cong
 \bigoplus_{\substack{\lambda\vdash n\\\ell(\lambda)\leq N}}
 \mathbb S_\lambda(U_N)\otimes(S^\lambda)^\vee.
\]
Tensoring with $\mathbb V^{\mathrm{tw}}_{g,k,n,d}(a;t)$ and taking the
$S_n$-invariant part yields
\begin{align*}
&\left[
 \mathbb V^{\mathrm{tw}}_{g,k,n,d}(a;t)
 \otimes U_N^{\otimes n}
\right]^{S_n}\\
&\cong
\bigoplus_{\substack{\lambda\vdash n\\\ell(\lambda)\leq N}}
\mathbb S_\lambda(U_N)\otimes
\operatorname{Hom}_{S_n}\left(
 S^\lambda,
 \mathbb V^{\mathrm{tw}}_{g,k,n,d}(a;t)
\right).
\end{align*}
Passing to virtual classes proves \eqref{mixed-Schur-Weyl-extraction}.  The
identity
\[
\left[
 \mathbb V^{\mathrm{tw}}_{g,k,n,d}(a;t)
 \otimes(S^\lambda)^\vee
\right]^{S_n}
\cong
\operatorname{Hom}_{S_n}\left(
 S^\lambda,
 \mathbb V^{\mathrm{tw}}_{g,k,n,d}(a;t)
\right)
\]
is precisely the invariant displayed in the theorem.  Finally, for $N\geq n$,
the classes
$[\mathbb S_\lambda(U_N)]$, $\lambda\vdash n$, are distinct irreducible
polynomial $GL_N$-characters and are therefore linearly independent.  Their
coefficients in the $GL_N$-equivariant invariant uniquely determine all
$m_\lambda$.
\end{proof}

For $N\geq n$, the $GL_N$-equivariant invariant therefore recovers the entire
virtual $S_n$-module.  Henceforth we simply call these quantities invariants;
the $GL_N$-equivariant interpretation is understood whenever its coefficients
are used.

The genus-$g$ mixed descendant potential is
\begin{equation*}
\mathcal F_g^{\mathrm{tw}}(x,t)
:=\sum_{k,n,d}\frac{Q^d}{k!}
\left\langle x(L),\dots,x(L);t(L),\dots,t(L)
\right\rangle^{S_n,\mathrm{tw}}_{g,n+k,d}.
\end{equation*}

\begin{mydef}
    The special case $N=1,x=0$, which extracts $S_n$-invariants of the sheaf cohomology of sheaves on $X_{g,n,d}$, is referred to as the \emph{symmetrized theory}. 
\end{mydef}

The distribution of permutable markings in boundary gluing formulas is
governed by the following equivariant binomial identity.

\begin{prop}[$S_n$-equivariant binomial identity]
\label{equivariant-binomial}
Let $n=n_1+\cdots+n_s$, put
$H=S_{n_1}\times\cdots\times S_{n_s}\subset S_n$, and let
$W_j$ be a virtual $S_{n_j}$-module.  For every $N\geq1$,
\begin{align}
\left[\operatorname{Ind}_{H}^{S_n}
(W_1\boxtimes\cdots\boxtimes W_s)\right]^{S_n}
&\cong W_1^{S_{n_1}}\otimes\cdots\otimes W_s^{S_{n_s}},
\label{equivariant-binomial-invariants}\\
\left[
\operatorname{Ind}_{H}^{S_n}(W_1\boxtimes\cdots\boxtimes W_s)
\otimes U_N^{\otimes n}
\right]^{S_n}
&\cong
\bigotimes_{j=1}^s
\left[W_j\otimes U_N^{\otimes n_j}\right]^{S_{n_j}}.
\label{equivariant-binomial-expanded}
\end{align}
The second identity is $GL_N$-equivariant, and both identities remain valid
with $K_T(\mathrm{pt})\otimes\Lambda$-coefficients.
\end{prop}

\begin{proof}
Frobenius reciprocity gives
\begin{align*}
\left[\operatorname{Ind}_{H}^{S_n}W\right]^{S_n}
&=\operatorname{Hom}_{S_n}
(\mathbf1,\operatorname{Ind}_{H}^{S_n}W)\\
&\cong\operatorname{Hom}_{H}(\mathbf1,W)=W^H,
\end{align*}
where $W=W_1\boxtimes\cdots\boxtimes W_s$.  Since the factors of $H$ act
on the corresponding tensor factors,
$W^H=\bigotimes_jW_j^{S_{n_j}}$, proving
\eqref{equivariant-binomial-invariants}.

For the second identity, apply the tensor identity for induction and then
Frobenius reciprocity:
\begin{align*}
&\left[
 \operatorname{Ind}_{H}^{S_n}W\otimes U_N^{\otimes n}
\right]^{S_n}\\
&\cong
\left[
 W\otimes\operatorname{Res}^{S_n}_{H}U_N^{\otimes n}
\right]^H.
\end{align*}
Under the ordered block decomposition of $\{1,\dots,n\}$ into consecutive
blocks of sizes $n_1,\dots,n_s$,
$\operatorname{Res}^{S_n}_{H}U_N^{\otimes n}
\cong U_N^{\otimes n_1}\boxtimes\cdots\boxtimes U_N^{\otimes n_s}$.
Consequently,
\begin{align*}
\left[
 W\otimes\operatorname{Res}^{S_n}_{H}U_N^{\otimes n}
\right]^H
&\cong
\left[
 \boxtimes_{j=1}^s(W_j\otimes U_N^{\otimes n_j})
\right]^{\prod_jS_{n_j}}\\
&\cong
\bigotimes_{j=1}^s
\left[W_j\otimes U_N^{\otimes n_j}\right]^{S_{n_j}},
\end{align*}
which proves \eqref{equivariant-binomial-expanded}.
\end{proof}

\subsection{Examples of type-I twistings}
The key examples of type-I twistings we consider are based on the
$K$-theoretic Euler class.

The $K$-theoretic Euler class is $\Eu(E):=\Lambda_{-1}(E^\vee)$.
For a line bundle $L$, $\Eu(L)=1-L^{-1}$ and, formally,
$\Eu(L)=\exp(\sum_{r<0}\Psi^r(L)/r)$.  We use
$\Eu_\lambda(L):=1-\lambda L^{-1}$ for the Euler class made invertible by an
auxiliary $\mathbb C^*$-character $\lambda$.  Multiplicativity and the
splitting principle define it for every virtual bundle.

A vector bundle $V$ on $X$ is \emph{convex} if
$H^1(C,f^*V)=0$ for every genus-zero stable map $f:C\to X$, and it is
\emph{concave} if $H^*(C,f^*V)=0$.  The target $X$ is convex when $T_X$ is
convex.  Convexity ensures that the relevant genus-zero index classes are
represented by honest bundles, although the virtual formulation does not
require this.

\subsubsection{Vector bundles}

Let $V\to X$ be a vector bundle and let $\mathbb C^*$ scale its fibers.  As
explained in Givental's discussion of Eulerian twistings
\cite[Corollary~2]{PermXI}, the fixed
locus in the stable-map moduli space of the total space of $V$ is
$X_{g,m,d}$, and its virtual normal class is $E_{g,m,d}(V)$.  Virtual
localization therefore gives
\begin{align*}
\chi^{\mathrm{vir}}_{\mathbb C^*}
\left(V_{g,m,d};\mathcal G\right)
={}&\chi^{\mathrm{vir}}_{\mathbb C^*}
\left(X_{g,m,d};
\mathcal G|_{X_{g,m,d}}\otimes
\Eu_\lambda(E_{g,m,d}(V))^{-1}\right).
\end{align*}
The same equality holds in the permutation-equivariant context because
localization commutes with permuting marked points.  Thus the equivariant
stable-map theory of the total space is the theory of $X$
twisted by the inverse Euler class.  In genus zero, when the relevant dual
bundle is convex, the nonequivariant limit exists.

\subsubsection{Complete intersections}

Let a regular section of a convex bundle $V\to X$ cut out a smooth complete
intersection $i:Y\hookrightarrow X$.   In genus 0, for convex $V$ the term $E_{0,m,d}(V)$, is an honest vector bundle, and the induced section of
$E_{0,m,d}(V)$ has the moduli space of maps to $Y$ as its zero locus.  Thus
$i_*\cO^{\mathrm{vir}}_{Y_{0,m,d}}
=\cO^{\mathrm{vir}}_{X_{0,m,d}}\otimes\Eu(E_{0,m,d}(V))$.
Before taking the limit $\lambda\to1$, we use $\Eu_\lambda$ so that the
twisting is invertible.  

\subsubsection{Level structures}

For an integer $\ell$ and a vector bundle $V$ on $X$, Ruan--Zhang define
level-$\ell$ quantum $K$-theory using the twisting
$\det(E_{g,m,d}(V))^{-\ell}$
\cite[Definition~2.4 and Remark~2.5]{ruan2018level}.

Over a general base algebra, the determinant cannot be written in the exponential form described previously, so to fit it into the theory, we use the following observation:

\begin{prop}\label{euler-determinant}
If $E$ has rank $r$ and its Euler classes are invertible, then
$\Eu(E)=(-1)^r\det(E)^{-1}\Eu(E^\vee)$.
\end{prop}

\begin{proof}
For a line bundle $L$ the ratio is
$(1-L^{-1})/(1-L)=-L^{-1}$.  If
$E=\bigoplus_{i=1}^rL_i$ after applying the splitting principle, then
$\Eu(E)/\Eu(E^\vee)=\prod_i(-L_i^{-1})
=(-1)^r\det(E)^{-1}$.
\end{proof}

Consequently, after adjoining two auxiliary characters $Y_+$ and $Y_-$, the
invertible class
\begin{equation*}
 (-1)^{\ell\operatorname{rk}E}
 \left(\frac{\Eu_{Y_+}(E)}{\Eu_{Y_-}(E^\vee)}\right)^\ell
\end{equation*}
specializes to $\det(E)^{-\ell}$ as $Y_+,Y_-\to1$. Degree-dependent signs
can equivalently be absorbed into the Novikov variables. 

\section{Basic Properties of Type I Theories}
\label{props}

Fix a type-I twisting $(C,V)$ and an integer $N\geq1$ as in
Definition~\ref{expanded-correlator-definition}.  In degree 0, the twisting
induces a pairing on $K_T(X)\otimes \Lambda$ given by $(a,b)^{\mathrm{tw}}:=\chi_T\bigl(X;a\otimes b\otimes C(V)\bigr).$

Let $\{\phi_\alpha\}$ be a basis of $K_T(X)\otimes\Lambda$, let
$g^{\mathrm{tw}}_{\alpha\beta}=(\phi_\alpha,\phi_\beta)^{\mathrm{tw}}$,
and write $g_{\mathrm{tw}}^{\alpha\beta}$ for the inverse matrix. Since $C$ is invertible, the pairing is nondegenerate. 

For named insertions $A_1,\dots,A_m$, define
\begin{equation*}
\llangle A_1,\dots,A_m\rrangle^{\mathrm{tw}}_{g,m;x,t}
:=\sum_{k,n\geq0}\sum_{d\in\operatorname{Eff}(X)}\frac{Q^d}{k!}
\left\langle A_1,\dots,A_m,x,\dots,x;t,\dots,t
\right\rangle^{S_n,\mathrm{tw}}_{g,m+k+n,d}.
\end{equation*}
There are exactly $k$ copies of the named input $x$ and $n$ copies
of the permutable input $t$ in each summand. Unless confusion is possible,
the subscripts $x,t$ will be suppressed.  Put
$G_{\alpha\beta}(x,t):=g^{\mathrm{tw}}_{\alpha\beta}
+\llangle\phi_\alpha,\phi_\beta\rrangle^{\mathrm{tw}}_{0,2;x,t}$, and denote
its inverse by $G^{\alpha\beta}(x,t)$.  The assumptions
$x,t\in\Lambda_+$ make the inverse well defined.

\subsection{Boundary gluing maps}

Let $M=I\sqcup J$ be a partition of the labeled set
$M=\{1,\dots,m\}$, with $g_1+g_2=g$ and $d_1+d_2=d$.
Write $\bullet$ and $\star$ for the two markings that will become the node.
The separating normalization space is
\begin{equation*}
Y_{I,d_1\mid J,d_2}
:=X_{g_1,I\cup\{\bullet\},d_1}
\mathop{\times}_{X}
X_{g_2,J\cup\{\star\},d_2}.
\end{equation*}
Here $X_{g_1,I\cup\{\bullet\},d_1}$ means the stable-map space with the
indicated labels and is isomorphic to
$X_{g_1,|I|+1,d_1}$.  The fiber product imposes
$\ev_\bullet=\ev_\star$.  Gluing the two domain curves at these markings
defines the nonequivariant boundary map
$\gl_{I,d_1\mid J,d_2}:Y_{I,d_1\mid J,d_2}\to X_{g,m,d}$.
We denote the common evaluation map on the fiber product by
$\ev_\Delta:Y_{I,d_1\mid J,d_2}\to X$.

For the nonseparating boundary, put
\begin{equation*}
Y_{\mathrm{irr}}
:=X_{g-1,M\cup\{\bullet,\star\},d}
\mathop{\times}_{X\times X,\,(\ev_\bullet,\ev_\star)}X,
\end{equation*}
where the map $X\to X\times X$ is the diagonal.  Thus
$Y_{\mathrm{irr}}=(\ev_\bullet,\ev_\star)^{-1}(\Delta_X)$, and gluing the
two branches gives
$\gl_{\mathrm{irr}}:Y_{\mathrm{irr}}\to X_{g,m,d}$.

We also need the corresponding maps for stabilized domain curves:
\begin{align*}
\Phi_{I\mid J}:&\
\overline M_{g_1,I\cup\{\bullet\}}
\times\overline M_{g_2,J\cup\{\star\}}
\longrightarrow\overline M_{g,M},\\
\Phi_{\mathrm{irr}}:&\
\overline M_{g-1,M\cup\{\bullet,\star\}}
\longrightarrow\overline M_{g,M}.
\end{align*}
The stabilization maps commute with $\gl_{I,d_1\mid J,d_2}$ and
$\gl_{\mathrm{irr}}$ on stable-map spaces and with $\Phi_{I\mid J}$ and
$\Phi_{\mathrm{irr}}$ on stabilized curves.

We now write the same maps with permutable markings.  Suppose the target
stable-map space is $X_{g,k+n,d}$, with named labels
$\{1,\dots,k\}=I\sqcup J$, and write $n=n_1+n_2$.  Set
\begin{align*}
Y^{n_1,n_2}_{I,d_1\mid J,d_2}:={}&
X_{g_1,|I|+n_1+1,d_1}\mathop{\times}_{X}
X_{g_2,|J|+n_2+1,d_2},\\
H_{n_1,n_2}:={}&S_{n_1}\times S_{n_2}\subset S_n.
\end{align*}
The first $S_{n_1}$ permutes the $n_1$ permutable markings on the first
factor and the second $S_{n_2}$ acts on the second factor.  The induced
$S_n$-space
\begin{equation*}
\widetilde Y^{n_1,n_2}_{I,d_1\mid J,d_2}
:=S_n\mathop{\times}_{H_{n_1,n_2}}
Y^{n_1,n_2}_{I,d_1\mid J,d_2}
\end{equation*}
records all choices of which $n_1$ of the $n$ labeled permutable markings lie
on the first component.  Gluing defines the $S_n$-equivariant map
$\gl^{S_n}_{I,d_1\mid J,d_2}:
\widetilde Y^{n_1,n_2}_{I,d_1\mid J,d_2}\to X_{g,k+n,d}$.
Passing to the $S_n$-quotients gives
\begin{equation*}
\left[Y^{n_1,n_2}_{I,d_1\mid J,d_2}/H_{n_1,n_2}\right]
\cong
\left[\widetilde Y^{n_1,n_2}_{I,d_1\mid J,d_2}/S_n\right]
\longrightarrow\left[X_{g,k+n,d}/S_n\right].
\end{equation*}

For genus reduction, define
$Y^{S_n}_{\mathrm{irr}}:=X_{g-1,k+n+2,d}
\mathop{\times}_{X\times X}X$,
where $S_n$ permutes the same $n$ markings and fixes the two branches of the
node.  The maps
$\gl^{S_n}_{\mathrm{irr}}:Y^{S_n}_{\mathrm{irr}}\to X_{g,k+n,d}$ and
$[Y^{S_n}_{\mathrm{irr}}/S_n]\to[X_{g,k+n,d}/S_n]$ are respectively the
$S_n$-equivariant gluing map and the map induced on the quotients.  This
distinguishes them from the nonequivariant maps $\gl_{I,d_1\mid J,d_2}$ and
$\gl_{\mathrm{irr}}$.

\subsection{Restriction of the twisting class}

The following restriction identities are due to Coates in
\cite[Lemmas~1.6.1--1.6.2]{Coatesthesis}. All pullbacks and pushforwards are taken $K$-theoretically. 

\begin{thm}[Restriction of the index class]\label{twist-restriction}
Let $E_{g,m,d}=R\ft_*\ev_{m+1}^*V$.
\begin{enumerate}[(i)]
\item If $\ft:X_{g,m+1,d}\to X_{g,m,d}$ forgets the last marked point, then
$E_{g,m+1,d}=\ft^*E_{g,m,d}$ and
$C(E_{g,m+1,d})=\ft^*C(E_{g,m,d})$; moreover,
$\cO^{\mathrm{vir},\mathrm{tw}}_{g,m+1,d}
=\ft^*\cO^{\mathrm{vir},\mathrm{tw}}_{g,m,d}$.
\item On the separating normalization space $Y_{I,d_1\mid J,d_2}$,
\begin{align}
\gl_{I,d_1\mid J,d_2}^*E_{g,m,d}
={}&\operatorname{pr}_1^*E_{g_1,|I|+1,d_1}
+\operatorname{pr}_2^*E_{g_2,|J|+1,d_2}
 -\ev_\Delta^*V,
\notag\\
\gl_{I,d_1\mid J,d_2}^*C(E_{g,m,d})
={}&\operatorname{pr}_1^*C(E_{g_1,|I|+1,d_1})
\otimes\operatorname{pr}_2^*C(E_{g_2,|J|+1,d_2})
\otimes\ev_\Delta^*C(V)^{-1}.
\label{twist-separating}
\end{align}
Thus the virtual gluing formula for the two twisted factors contains the node
factor $\ev_\Delta^*C(V)^{-1}$.
\item On $Y_{\mathrm{irr}}$,
\begin{align*}
\gl_{\mathrm{irr}}^*E_{g,m,d}
&=E_{g-1,m+2,d}-\ev_\Delta^*V,\\
\gl_{\mathrm{irr}}^*C(E_{g,m,d})
&=C(E_{g-1,m+2,d})\otimes\ev_\Delta^*C(V)^{-1}.
\end{align*}
\end{enumerate}
These identities are $S_n$-equivariant on
$\widetilde Y^{n_1,n_2}_{I,d_1\mid J,d_2}$ and
$Y^{S_n}_{\mathrm{irr}}$ and therefore descend through the displayed quotient
maps without alteration.
\end{thm}

\subsection{Kontsevich--Manin axioms}

Define
\begin{equation*}
H^{\mathrm{tw}}_{g,m,d}(A_1,\dots,A_m)
:=\cO^{\mathrm{vir},\mathrm{tw}}_{g,m,d}\otimes
\prod_{i=1}^m\ev_i^*A_i(L_i).
\end{equation*}

For a stable pair $(g,m)$, let
$\st_{k,n,d}:X_{g,m+k+n,d}\to\overline M_{g,m}$ forget the map and the last
$k+n$ markings and stabilize the domain.  If
$\mathcal E$ is $S_n$-equivariant, define
\begin{align*}
(\st_{k,n,d})_*^{S_n}\mathcal E
&:=\left[
R(\st_{k,n,d})_*(\mathcal E\otimes U_N^{\otimes n})
\right]^{S_n}
\\
&\in K_{T\times GL_N}(\overline M_{g,m})\otimes\Lambda.
\end{align*}
The $S_n$-action on $U_N^{\otimes n}$ is the one in
Definition~\ref{expanded-correlator-definition}.  Define
\begin{align*}
\mathbb I^{\mathrm{tw}}_{g,m}
(A_1,\dots,A_m;x,t)
:={}&\sum_{k,n\geq0}\sum_{d\in\operatorname{Eff}(X)}
\frac{Q^d}{k!}(\st_{k,n,d})_*^{S_n}
\Bigl(
H^{\mathrm{tw}}_{g,m+k+n,d}
(A_1,\dots,A_m,\underbrace{x,\dots,x}_{k};
\underbrace{t,\dots,t}_{n})\Bigr).
\end{align*}
Here the semicolon only separates the $n$ permutable insertions from the
$m+k$ named insertions.  Applying the equivariant Euler
characteristic gives exactly the invariants described by double brackets:
\begin{equation*}
\chi_{T\times GL_N}\!\left(\overline M_{g,m};
\mathbb I^{\mathrm{tw}}_{g,m}(A_1,\dots,A_m;x,t)\right)
=\llangle A_1,\dots,A_m\rrangle^{\mathrm{tw}}_{g,m;x,t}.
\end{equation*}

The following are the Kontsevich-Manin axioms for twisted $S_n$-equivariant $K$-theoretic invariants: 

\begin{thm}\label{mixed-axioms}
\hfill
\begin{enumerate}[(i)]
\item \emph{Covariance.} For every $\sigma\in S_m$,
$\sigma^*H^{\mathrm{tw}}_{g,m,d}(A_1,\dots,A_m)
=H^{\mathrm{tw}}_{g,m,d}
(A_{\sigma(1)},\dots,A_{\sigma(m)})$.
In particular, equal insertions at the last $n$ markings give an
$S_n$-equivariant class.
\item \emph{Fundamental class.} If
$\ft:X_{g,m+1,d}\to X_{g,m,d}$ forgets the last named point and the
$A_i$ are primary, then
$H^{\mathrm{tw}}_{g,m+1,d}(A_1,\dots,A_m,1)
=\ft^*H^{\mathrm{tw}}_{g,m,d}(A_1,\dots,A_m)$.
\item \emph{Mapping to a point.} Let $\mathcal H$ be the Hodge bundle, so
$\mathcal H^\vee=R^1(\pi_{\mathcal C})_*\cO_{\mathcal C}$.  Then
\begin{align}
 X_{g,m,0}&=\overline M_{g,m}\times X,\notag\\
 \cO^{\mathrm{vir},\mathrm{tw}}_{g,m,0}
 &=\lambda_{-1}\!\left(
 \mathcal H^\vee\boxtimes T_X\right)^\vee
 \otimes C\!\left((1-\mathcal H^\vee)\boxtimes V\right).
 \label{mapping-point-twisted-sheaf}
\end{align}
In particular,
$\langle a,b,c\rangle^{\mathrm{tw}}_{0,3,0}=\chi_T(X;abcC(V))$.
\item \emph{Splitting.} For $M=I\sqcup J$ and the gluing map
$\Phi_{I\mid J}$,
\begin{align}
\Phi_{I\mid J}^*\mathbb I^{\mathrm{tw}}_{g,m}(A_M;x,t)
={}&\sum_{\alpha,\beta}
\mathbb I^{\mathrm{tw}}_{g_1,|I|+1}(A_I,\phi_\alpha;x,t)
\boxtimes
\mathbb I^{\mathrm{tw}}_{g_2,|J|+1}(\phi_\beta,A_J;x,t)
G^{\alpha\beta}(x,t).
\label{mixed-splitting}
\end{align}

\item \emph{Genus reduction.} For $\Phi_{\mathrm{irr}}$,
\begin{align*}
 \Phi_{\mathrm{irr}}^*\mathbb I^{\mathrm{tw}}_{g,m}(A_M;x,t)
={}&\sum_{\alpha,\beta}
\mathbb I^{\mathrm{tw}}_{g-1,m+2}
(A_M,\phi_\alpha,\phi_\beta;x,t)
G^{\alpha\beta}(x,t).
\end{align*}
\end{enumerate}
\end{thm}

The covariance, fundamental-class, and mapping-to-a-point assertions follow
from the untwisted arguments of \cite{yplee}, using
$E_{g,m+1,d}=\ft^*E_{g,m,d}$ and
$E_{g,m,0}=(1-\mathcal H^\vee)\boxtimes V$.

\subsection{Splitting and genus reduction}\label{split}

Fix $M=I\sqcup J$.  For a summand in the definition of
$\mathbb I^{\mathrm{tw}}_{g,m}$, the inverse image of the boundary
$\operatorname{im}(\Phi_{I\mid J})$ under $\st_{k,n,d}$ is a virtual
normal-crossings divisor.  A depth-$r$ normalization stratum has two
components carrying the labels in $I$ and $J$, joined by a chain with $r$
nodes and $r-1$ intermediate genus-zero components.  Choose data satisfying
$\sum_{a=0}^{r}k_a=k$, $\sum_{a=0}^{r}n_a=n$, and
$\sum_{a=0}^{r}d_a=d$,
and set
\begin{align*}
Y_r(\mathbf k,\mathbf n,\mathbf d):={}&
X_{g_1,|I|+k_0+n_0+1,d_0}
\mathop{\times}_X
\left(\mathop{\times}_{a=1}^{r-1}{}_X
X_{0,k_a+n_a+2,d_a}\right)
\mathop{\times}_X
X_{g_2,|J|+k_r+n_r+1,d_r}.
\end{align*}
Every fiber product defining $Y_r(\mathbf k,\mathbf n,\mathbf d)$ identifies the two
evaluations at consecutive branches of a node.  The usual stability
conditions are imposed on every factor.  Put
$H_{\mathbf n}=\prod_{a=0}^rS_{n_a}$.  The $S_n$-equivariant depth-$r$ stratum
is $S_n\mathop{\times}_{H_{\mathbf n}}
Y_r(\mathbf k,\mathbf n,\mathbf d)$, and its gluing map to
$X_{g,m+k+n,d}$ is $S_n$-equivariant.  This induced space is the iterated
version of $\widetilde Y^{n_1,n_2}_{I,d_1\mid J,d_2}$ and records all
distributions of the $n$ permutable labels among the $r+1$ components.

If $\iota_r:Y_r(\mathbf k,\mathbf n,\mathbf d)\to X_{g,m+k+n,d}$ is the
nonequivariant normalization map, Theorem~\ref{twist-restriction}, applied at
each of its $r$ nodes, gives
\begin{equation*}
 \iota_r^*C(E_{g,m+k+n,d})
 =\left(\boxtimes_{a=0}^{r}C(E_a)\right)
 \otimes\left(\bigotimes_{e=1}^{r}\ev_e^*C(V)^{-1}\right).
\end{equation*}
For one node the K\"unneth expansion of the resulting kernel is
\begin{equation}\label{twisted-node-kernel}
 \Delta_*C(V)^{-1}
 =\sum_{\alpha,\beta}g_{\mathrm{tw}}^{\alpha\beta}
 \phi_\alpha\boxtimes\phi_\beta.
\end{equation}
Indeed, the right-hand side is the inverse tensor for
$g^{\mathrm{tw}}_{\alpha\beta}
=\chi_T(X;\phi_\alpha\phi_\beta C(V))$; the factor $C(V)^{-1}$ is precisely
the last factor of \eqref{twist-separating}.

Let $D_r$ denote the union over the preceding distribution data of the
depth-$r$ induced strata.  The formula for the structure sheaf of a virtual
normal-crossings divisor gives
\begin{equation}\label{boundary-inclusion-exclusion}
 [\cO_{\st^{-1}\operatorname{im}(\Phi_{I\mid J})}^{\mathrm{vir},\mathrm{tw}}]
 =\sum_{r\geq1}(-1)^{r-1}[\cO_{D_r}^{\mathrm{vir},\mathrm{tw}}].
\end{equation}
Put
\begin{equation*}
 F_{\alpha\beta}
 :=\llangle\phi_\alpha,\phi_\beta\rrangle^{\mathrm{tw}}_{0,2;x,t}
 =\sum_{k,n\geq0}\sum_d\frac{Q^d}{k!}
 \left\langle\phi_\alpha,\phi_\beta,
 \underbrace{x,\dots,x}_{k};
\underbrace{t,\dots,t}_{n}\right\rangle^{S_n,\mathrm{tw}}_{0,k+n+2,d}.
\end{equation*}
Equations \eqref{twisted-node-kernel} and
\eqref{boundary-inclusion-exclusion} show that a chain with $r$ nodes
contributes
$(-1)^{r-1}g_{\mathrm{tw}}^{-1}
(Fg_{\mathrm{tw}}^{-1})^{r-1}$.  Therefore the sum over all chain lengths is
\begin{align*}
 &g_{\mathrm{tw}}^{-1}
 -g_{\mathrm{tw}}^{-1}Fg_{\mathrm{tw}}^{-1}
 +g_{\mathrm{tw}}^{-1}Fg_{\mathrm{tw}}^{-1}Fg_{\mathrm{tw}}^{-1}-\cdots
 \\
 &\hspace{35mm}=(g^{\mathrm{tw}}+F)^{-1}=G^{-1}.
\end{align*}
This calculation takes place in
$K_T(\mathrm{pt})\otimes\Lambda\otimes R(GL_N)_{\mathbb Q}$.

For the $k$ additional named markings, the distribution gives
$k!^{-1}\binom{k}{k_0,\dots,k_r}=\prod_{a=0}^{r}k_a!^{-1}$.
For the permutable markings, let $W_{n_a}$ be the virtual cohomology module
of the $a$th component after all its insertions and node classes have been
included, but before taking $S_{n_a}$-invariants.  Because the depth-$r$
stratum is induced from $H_{\mathbf n}=\prod_aS_{n_a}$, the virtual module
contributing to that stratum is
$\operatorname{Ind}_{H_{\mathbf n}}^{S_n}
(\boxtimes_aW_{n_a})$.  After inserting $U_N$ at every permutable marking,
Definition~\ref{expanded-correlator-definition} identifies its contribution
to the correlator with the left-hand side of
\begin{equation}
\left[
\operatorname{Ind}_{H_{\mathbf n}}^{S_n}
\left(\boxtimes_{a=0}^{r}W_{n_a}\right)
\otimes U_N^{\otimes n}
\right]^{S_n}
\cong
\bigotimes_{a=0}^{r}
\left[W_{n_a}\otimes U_N^{\otimes n_a}\right]^{S_{n_a}}.
\label{expanded-depth-r-equivariant-binomial}
\end{equation}
By Proposition~\ref{equivariant-binomial}, the two sides are equal.  Each
factor on the right is exactly the contribution of the
corresponding component with $n_a$ permutable markings. 
Substitution of the preceding
multinomial identity and \eqref{expanded-depth-r-equivariant-binomial} into
\eqref{boundary-inclusion-exclusion} yields
\begin{align*}
\Phi_{I\mid J}^*\mathbb I^{\mathrm{tw}}_{g,m}(A_M;x,t)
={}&\sum_{\alpha,\beta}
\mathbb I^{\mathrm{tw}}_{g_1,|I|+1}(A_I,\phi_\alpha;x,t)
\boxtimes
\mathbb I^{\mathrm{tw}}_{g_2,|J|+1}(\phi_\beta,A_J;x,t)
\notag\\[-2mm]
&\quad\times
\left(\sum_{r\geq1}(-1)^{r-1}
g_{\mathrm{tw}}^{-1}(Fg_{\mathrm{tw}}^{-1})^{r-1}
\right)^{\alpha\beta}
\notag\\
={}&\sum_{\alpha,\beta}
\mathbb I^{\mathrm{tw}}_{g_1,|I|+1}(A_I,\phi_\alpha;x,t)
\boxtimes
\mathbb I^{\mathrm{tw}}_{g_2,|J|+1}(\phi_\beta,A_J;x,t)
G^{\alpha\beta},
\end{align*}
which proves \eqref{mixed-splitting}.

\subsection{String, dilaton, and WDVV equations}

Let $\ft:X_{g,m+1,d}\to X_{g,m,d}$ forget the last named marked
point, and let $\mathcal H$ be the Hodge bundle.  The following identities are
identities of $S_n$-equivariant $K$-classes whenever $S_n$ permutes a disjoint
set of markings.  Their untwisted forms are 
\cite[\S\S4.4--4.5]{yplee}; the  permutation-equivariant potential
identities are 
\cite[Propositions~1--2]{PermVII}.

\begin{thm}\label{string-dilaton}
For genus zero,
\begin{equation}\label{sheaf-string-zero}
\ft_*\left(\cO^{\mathrm{vir},\mathrm{tw}}
\prod_{i=1}^{m}\frac{1}{1-q_iL_i}\right)
=\cO^{\mathrm{vir},\mathrm{tw}}
\left(1+\sum_{i=1}^{m}\frac{q_i}{1-q_i}\right)
\prod_{i=1}^{m}\frac{1}{1-q_iL_i}.
\end{equation}
For $g\geq1$,
\begin{align*}
&\ft_*\left(\cO^{\mathrm{vir},\mathrm{tw}}
\frac{1}{1-q\mathcal H^\vee}
\prod_{i=1}^{m}\frac{1}{1-q_iL_i}\right)\notag\\
&=\cO^{\mathrm{vir},\mathrm{tw}}
\frac{1}{1-q\mathcal H^\vee}
\left(1-\mathcal H^\vee+
\sum_{i=1}^{m}\frac{q_i}{1-q_i}\right)
\prod_{i=1}^{m}\frac{1}{1-q_iL_i}.
\end{align*}
The dilaton identity is
\begin{align*}
&\ft_*\left(\cO^{\mathrm{vir},\mathrm{tw}}
\frac{1}{1-q\mathcal H}
\prod_{i=1}^{m}\frac{1}{1-q_iL_i}\,L_{m+1}\right)\notag\\
&=\cO^{\mathrm{vir},\mathrm{tw}}
\frac{1}{1-q\mathcal H}
\left(\mathcal H-1+
\sum_{i=1}^{m}\frac{1}{1-q_i}\right)
\prod_{i=1}^{m}\frac{1}{1-q_iL_i}.
\end{align*}
In particular,
$\ft_*(L_{m+1}-1)=\mathcal H+\mathcal H^\vee+(m-2)\mathbf1$.
\end{thm}

\begin{proof}
 For any class $F$, Theorem~\ref{twist-restriction}(i)
and the projection formula give
\begin{align*}
 \ft_*\bigl(C(E_{g,m+1,d})\otimes F\bigr)
 &=\ft_*\bigl(\ft^*C(E_{g,m,d})\otimes F\bigr)\\
 &=C(E_{g,m,d})\otimes\ft_*F.
\end{align*}
Thus tensoring the two sides of 
\cite[\S\S4.4--4.5]{yplee} by $C(E_{g,m,d})$ gives the three displayed
identities.  If $\sigma\in S_n$ acts on markings disjoint from the forgotten
point, then $\ft\circ\sigma=\sigma\circ\ft$ and
$\sigma^*C(E_{g,m+1,d})=C(E_{g,m+1,d})$,
so these equalities already hold in $S_n$-equivariant $K$-theory.
\end{proof}

For genus zero, \eqref{sheaf-string-zero} is equivalently
\begin{align}\label{descendant-string}
&\left\langle a_1(L),\dots,a_m(L),1;t,\dots,t
\right\rangle^{S_n,\mathrm{tw}}_{0,m+n+1,d}\notag\\
&\quad=\left\langle a_1(L),\dots,a_m(L);t,\dots,t
\right\rangle^{S_n,\mathrm{tw}}_{0,m+n,d}\notag\\
&\quad+
\sum_{i=1}^m\left\langle a_1(L),\dots,
\frac{a_i(L)-a_i(1)}{L-1},\dots,a_m(L);t,\dots,t
\right\rangle^{S_n,\mathrm{tw}}_{0,m+n,d}.
\end{align}
For primary $a_i$ the sum vanishes, giving the fundamental-class identity.


\begin{thm}\label{-wdvv}
Denote the twisted quantum metric by $G$.
For $a,b\in K_T(X)\otimes\Lambda$ and formal variables $u,v$,
\begin{align}
&(a,b)^{\mathrm{tw}}
+(1-uv)\llangle\frac{a}{1-uL},\frac{b}{1-vL}
\rrangle^{\mathrm{tw}}_{0,2}\notag\\
&=\sum_{\alpha,\beta}
\left((a,\phi_\alpha)^{\mathrm{tw}}
+\llangle\frac{a}{1-uL},\phi_\alpha\rrangle^{\mathrm{tw}}_{0,2}\right)
G^{\alpha\beta}
\left((\phi_\beta,b)^{\mathrm{tw}}
+\llangle\phi_\beta,\frac{b}{1-vL}\rrangle^{\mathrm{tw}}_{0,2}\right).
\label{-wdvv-equation}
\end{align}
\end{thm}

\begin{proof}
Let $\mathrm{st}:X_{0,4+k+n,d}\to\overline M_{0,4}\cong\mathbb P^1$ retain
the four named markings.  For the boundary points
$0,1,\infty\in\overline M_{0,4}$, we have
$[\cO_0]=[\cO_1]=[\cO_\infty]$ in $K(\overline M_{0,4})$,
and hence
\[
 \bigl[\cO^{\mathrm{vir},\mathrm{tw}}\otimes^{\mathbf L}
 \mathrm{st}^*\cO_0\bigr]
 =
 \bigl[\cO^{\mathrm{vir},\mathrm{tw}}\otimes^{\mathbf L}
 \mathrm{st}^*\cO_1\bigr]
 =
 \bigl[\cO^{\mathrm{vir},\mathrm{tw}}\otimes^{\mathbf L}
 \mathrm{st}^*\cO_\infty\bigr].
\]
Using \eqref{mixed-splitting} on two of these fibers gives
\begin{align*}
 &\sum_{\mu,\nu}
 \llangle a,b,\phi_\mu\rrangle^{\mathrm{tw}}
 G^{\mu\nu}
 \llangle\phi_\nu,c,d\rrangle^{\mathrm{tw}}\\
 &=
 \sum_{\mu,\nu}
 \llangle a,c,\phi_\mu\rrangle^{\mathrm{tw}}
 G^{\mu\nu}
 \llangle\phi_\nu,b,d\rrangle^{\mathrm{tw}}.
\end{align*}

For a decomposition $n=n_1+n_2$ of the permutable markings,
the virtual module on the normalization of the corresponding boundary
divisor is
$\operatorname{Ind}_{S_{n_1}\times S_{n_2}}^{S_n}(W_1\boxtimes W_2)$.
After tensoring with $U_N^{\otimes n}$ and taking $S_n$-invariants,
\eqref{equivariant-binomial-expanded} identifies the contribution of this
boundary stratum to the invariant with
$[W_1\otimes U_N^{\otimes n_1}]^{S_{n_1}}\otimes
[W_2\otimes U_N^{\otimes n_2}]^{S_{n_2}}$.  These two factors are exactly
the two correlators appearing on that side of WDVV, so summing over
$n_1+n_2=n$ produces the product of double brackets without an additional
binomial coefficient.
Substituting
the two descendant insertions and applying \eqref{descendant-string} gives
\eqref{-wdvv-equation}.  This is precisely Givental's proof of
\cite[Proposition~3]{PermVII}, with the twisting factor included.
\end{proof}

\section{Twisted Quantum Products}
\label{product}

Fix primary inputs $x,t\in\Lambda_+$, with $t$ permutable and
$x=\sum_\alpha x^\alpha\phi_\alpha$ named.  The
\emph{ twisted quantum product} $\star_{x,t}$ is defined by
\begin{equation*}
 \bigl(\phi_\alpha\star_{x,t}\phi_\beta,\phi_\gamma
 \bigr)_{x,t}
 :=\llangle\phi_\alpha,\phi_\beta,\phi_\gamma
 \rrangle^{\mathrm{tw}}_{0,3;x,t},
\end{equation*}
where $(\phi_\alpha,\phi_\beta)_{x,t}=G_{\alpha\beta}(x,t)$.  Thus
\begin{equation*}
 \phi_\alpha\star_{x,t}\phi_\beta
 =\sum_{\mu,\nu}
\llangle\phi_\alpha,\phi_\beta,\phi_\mu
 \rrangle^{\mathrm{tw}}_{0,3;x,t}
 G^{\mu\nu}(x,t)\phi_\nu.
\end{equation*}

\begin{thm}\label{-frobenius-ring}
The product $\star_{x,t}$ is commutative and associative with unit 1. 
\end{thm}

\begin{proof}
Associativity is equivalent to
\begin{align*}
 &\sum_{\mu,\nu}
 \llangle a,b,\phi_\mu\rrangle^{\mathrm{tw}}
 G^{\mu\nu}
 \llangle\phi_\nu,c,d\rrangle^{\mathrm{tw}}\\
 &=
 \sum_{\mu,\nu}
 \llangle a,c,\phi_\mu\rrangle^{\mathrm{tw}}
 G^{\mu\nu}
 \llangle\phi_\nu,b,d\rrangle^{\mathrm{tw}},
\end{align*}
which is the four-point consequence of Theorem~\ref{-wdvv}.

For the unit, \eqref{descendant-string} and
\eqref{mapping-point-twisted-sheaf} give
\begin{align*}
 \llangle1,a,b\rrangle^{\mathrm{tw}}
 &=\langle1,a,b\rangle^{\mathrm{tw}}_{0,3,0}
   +\llangle a,b\rrangle^{\mathrm{tw}}_{0,2}\\
 &=\chi_T(X;abC(V))
   +\llangle a,b\rrangle^{\mathrm{tw}}_{0,2}
 =G(a,b).
\end{align*}
Therefore $G(1\star_{x,t}a,b)=G(a,b)$ and $1\star_{x,t}a=a$.
\end{proof}

\section{$q$-difference Modules}
\label{qdiff}

All constructions in this section are made at a fixed permutable input $t$, whereas we let the nonpermuted input $x$ vary. We establish some results about generating functions of twisted invariants and their relationship to quantum $K$-rings. These results are all originally due to Iritani-Milanov-Tonita in \cite{imt} and adapted to the twisted and $S_n$-equivariant case. 


\subsection{Quantum Connection}

Write $\partial_\alpha=\partial/\partial x^\alpha$.  The  quantum
connection is $\nabla_\alpha:=(1-q)\partial_\alpha+
\phi_\alpha\star_{x,t}$.
Associativity implies $[\nabla_\alpha,\nabla_\beta]=0$.  Define the $S$-matrix
$S_{x,t}(q)$ by
\begin{equation*}
 \bigl(S_{x,t}(q)a,b\bigr)^{\mathrm{tw}}
 =(a,b)^{\mathrm{tw}}+
 \llangle\frac{a}{1-qL},b\rrangle^{\mathrm{tw}}_{0,2;x,t}.
\end{equation*}
It is invertible in the $(Q,x,t)$-adic topology.  With
$T_{x,t}=S_{x,t}^{-1}$, the WDVV and string equations give
$(1-q)\partial_\alpha T_{x,t}=T_{x,t}(\phi_\alpha\star_{x,t})$ and
$T_{x,t}=1+O(Q,x,t)$, as in \cite[Theorem~2.2]{imt}.  Thus $T_{x,t}$ is a
fundamental solution of the quantum connection.  The proof of
\cite[Proposition~2.3]{imt} also gives
$S_{x,t}(q)^*S_{x,t}(q^{-1})=1$.

\subsection{ $J$-functions and loop spaces}

Let $\Kloop^{\mathrm{tw}}$ be the space of rational functions of $q$ with
coefficients in $K_T(X)\otimes\Lambda$, completed as above, and equipped
with the symplectic form
\begin{equation*}
 \Omega^{\mathrm{tw}}(f,g)
 =-\operatorname*{Res}_{q=0,\infty}
 (f(q),g(q^{-1}))^{\mathrm{tw}}\frac{dq}{q}.
\end{equation*}
It has a Lagrangian polarization
$\Kloop^{\mathrm{tw}}=\Kloop_+\oplus\Kloop_-$, where $\Kloop_+$ consists of
Laurent polynomials and $\Kloop_-$ consists of rational functions regular at
$q=0$ and vanishing at $q=\infty$. 

For descendant inputs $x(q),t(q)\in\Kloop_+$, the  twisted big
$J$-function is
\begin{align*}
 \J^{\mathrm{tw}}(x,t)
 :={}&1-q+x(q)+t(q)\notag\\
 &+\sum_{\alpha,k,n,d}\frac{Q^d}{k!}\,\phi_\alpha
 \left\langle
 \frac{\phi^\alpha}{1-qL},x(L),\dots,x(L);
 t(L),\dots,t(L)
 \right\rangle^{S_n,\mathrm{tw}}_{0,1+k+n,d}.
\end{align*}
Here $\{\phi^\alpha\}$ is dual to $\{\phi_\alpha\}$ for the twisted
pairing, and only stable correlators are included in the sum.  

We use Givental's notation for two related ranges.  Define
$\mathcal L^{\mathrm{tw}}$ to be the range of
$t\mapsto\J^{\mathrm{tw}}(0,t)$ for $t\in\Kloop_+$ with
$t(1)\in\Lambda_+$.  Thus the named input is $x=0$, while the
permutable input $t$ varies.  For a fixed permutable input $t$, define
$\mathcal L_t^{\mathrm{tw}}$ to be the range of
$x\mapsto\J^{\mathrm{tw}}(x,t)$.  Thus the subscript $t$ always means that
$t$ is held fixed while $x$ varies.  We use the same convention below after
adding labels for the target and the twisting.

By Kawasaki--Riemann--Roch, the output of every correlator lies in
$\Kloop_-$, so projection to $\Kloop_+$ recovers the input.

For primary inputs, $\J^{\mathrm{tw}}(x,t)=(1-q)T_{x,t}(q)1$.

\begin{thm}
\label{-cone-theorem}
For each fixed primary permutable input $t\in\Lambda_+$,
$\mathcal L_t^{\mathrm{tw}}\subset\Kloop^{\mathrm{tw}}$ is a formal
overruled Lagrangian cone.  If $T_f\mathcal L_t^{\mathrm{tw}}$ is its
tangent space at $f$, then
\[
 (1-q)T_f\mathcal L_t^{\mathrm{tw}}\subset\mathcal L_t^{\mathrm{tw}},
 \quad\text{and}\quad
 T_g\mathcal L_t^{\mathrm{tw}}=T_f\mathcal L_t^{\mathrm{tw}}
 \quad(g\in(1-q)T_f\mathcal L_t^{\mathrm{tw}}).
\]
The tangent space at $\J^{\mathrm{tw}}(x,t)$ is
$T_{x,t}(\Kloop_+)$.  The unsubscripted range
$\mathcal L^{\mathrm{tw}}$, obtained by setting $x=0$ and varying $t$, is
overruled in the permutation-equivariant sense, but it is not in general
Lagrangian.
\end{thm}

\begin{proof}
    This is identical to the untwisted case proved by Givental in \cite[Theorem~1]{PermVII}.
\end{proof}

Recall that for torus equivariant parameters $\Lambda_i$, coming from  a torus acting on $X$, we complete at
$1-\Lambda_i$.  However parameters associated with $\mathbb{C}^*$-actions on vector bundles, e.g. $\lambda$ in the definition of $Eu_\lambda$, are not completed, instead, they are handled in Yan's rational
loop space: coefficients are rational in these parameters and the positive
and negative spaces are separated by their allowed poles, as in
\cite[\S3]{Xiaohanserre}.  This avoids incompatible simultaneous completions and
does not change the arguments below.

\subsection{$q$-shift operators}

Choose line bundles $P_1,\dots,P_r$ and Novikov variables
$Q_1,\dots,Q_r$ so that
$-c_1(P_i)$ is dual to $Q_i$.  The ruling spaces of
$\mathcal L_t^{\mathrm{tw}}$ are modules over the finite-difference algebra
generated by $P_iq^{Q_i\partial_{Q_i}}$ and $Q^d$,
and by multiplication with Laurent polynomials in $q$, by the argument as \cite[Theorem~6.6]{tonitaktwist}.  

Conjugating by the $S$-matrix defines the shift operator $\mathcal{A}_i$:
\begin{equation*}
 \mathcal A_i
 :=S_{x,t}\circ P_iq^{Q_i\partial_{Q_i}}\circ S_{x,t}^{-1}
 =T_{x,t}^{-1}\circ P_iq^{Q_i\partial_{Q_i}}\circ T_{x,t}
 =A_i(q,Q,x,t)q^{Q_i\partial_{Q_i}}.
\end{equation*}
Here $A_i(q,Q,x,t)$ is an endomorphism-valued coefficient; it is not a
$K$-theory class.

These operators satisfy the following properties:
\begin{thm}
The coefficient $A_i$ is polynomial in $q$ after the chosen completion and
satisfies $[\mathcal A_i,\mathcal A_j]=0$ and
$[\nabla_\alpha,\mathcal A_i]=0$.
\end{thm}

The proof is identical to the proofs of \cite[Prop~2.6--2.10]{imt}. 

 Define 
 $A_{i,\mathrm{com}}(x,t)
 :=\left.\mathcal A_i\right|_{q=1}
 =A_i(1,Q,x,t)
 \in\operatorname{End}_{\Lambda}
 \bigl(K_T(X)\otimes\Lambda\bigr).$ We then have:

\begin{cor}\label{commuting-shift-is-multiplication}
The endomorphism $A_{i,\mathrm{com}}(x,t)$ commutes with every quantum
multiplication operator.  Define the quantum line-bundle \emph{class} by
$\widehat P_i(x,t):=A_{i,\mathrm{com}}(x,t)1\in K_T(X)\otimes\Lambda$.
Then, for every $v\in K_T(X)\otimes\Lambda$,
\begin{equation}\label{commuting-shift-multiplication}
 A_{i,\mathrm{com}}(x,t)v
 =v\star_{x,t}\widehat P_i(x,t),
 \quad
 A_{i,\mathrm{com}}(x,t)
 =\widehat P_i(x,t)\star_{x,t},
\end{equation}
and $\widehat P_i=P_i+O(Q)$.
\end{cor}

\begin{proof}
Specializing $[\nabla_\alpha,\mathcal A_i]=0$ at $q=1$ gives
$[\phi_\alpha\star_{x,t},A_{i,\mathrm{com}}]=0$; this is
\cite[Cor~2.9]{imt}.  Hence
$A_{i,\mathrm{com}}v
=A_{i,\mathrm{com}}((v\star_{x,t})1)
=(v\star_{x,t})A_{i,\mathrm{com}}1
=v\star_{x,t}\widehat P_i$.
Thus $A_{i,\mathrm{com}}$ acts by quantum multiplication by
$A_{i,\mathrm{com}}1=P_i+O(Q)$.
\end{proof}

If $F$ is a polynomial and $c_1,\dots,c_s$ are classes, the unadorned
expression $F(c_1,\dots,c_s)$ denotes evaluation in the classical $K$-ring,
using ordinary tensor-product multiplication.  For evaluation using the
quantum product, write
$F^{\star_{x,t}}(c_1,\dots,c_s):=
F(c_1\star_{x,t},\dots,c_s\star_{x,t})1$.
Thus
$F(A_{1,\mathrm{com}},\dots,A_{s,\mathrm{com}})1
=F^{\star_{x,t}}(\widehat P_1,\dots,\widehat P_s)$; the left-hand side
evaluates a polynomial in endomorphisms, while the right-hand side is a class
in the quantum ring.

It is sometimes convenient to remove $K$-theory classes from the difference
operators.  Put
\begin{equation*}
 \widetilde\J^{\mathrm{tw}}
 :=\left(\prod_iP_i^{\log Q_i/\log q}\right)\J^{\mathrm{tw}},
 \quad
 \widetilde S:=S\left(\prod_iP_i^{-\log Q_i/\log q}\right).
\end{equation*}
The powers are defined by the finite binomial expansion in the nilpotent
classes $1-P_i$.  Conjugating $q^{Q_i\partial_{Q_i}}$ by the logarithmic
factor contributes multiplication by $P_i$, so we have
\begin{equation}\label{normalized-intertwining}
 q^{Q_i\partial_{Q_i}}\widetilde S^{-1}
 =\widetilde S^{-1}\mathcal A_i,
 \quad
 (1-q)\partial_\alpha\widetilde S^{-1}
 =\widetilde S^{-1}(\phi_\alpha\star_{x,t}).
\end{equation}
\subsection{Difference Equations}

Let $D$ be a polynomial operator, regular at $q=1$, in $q,Q$,  the shifts
$q^{Q_i\partial_{Q_i}}$, and the operators
$(1-q)\partial_\alpha$ with scalar coefficients.  In the following
statement, replacing an operator by an endomorphism means preserving the
order of the factors.

The following theorem is the type-I twisted,  extension of the
untwisted non-equivariant relation theorem
\cite[Proposition~2.12]{imt}.

\begin{thm}\label{rel}
If $D(q,Q,q^{Q_i\partial_{Q_i}},(1-q)\partial_\alpha)
\widetilde\J^{\mathrm{tw}}(x,t)=0$,
then the  twisted quantum ring satisfies
\begin{equation}\label{relation-theorem-equation}
 D\bigl(1,Q,A_{i,\mathrm{com}}(x,t),
 \phi_\alpha\star_{x,t}\bigr)1
 =D^{\star_{x,t}}\bigl(1,Q,\widehat P_i(x,t),\phi_\alpha\bigr)=0.
\end{equation}
The same conclusion holds if we remove dependence on $\phi_\alpha$, and fix the value of $x$. 
Equivalently, without logarithmic normalization,
\begin{align*}
 &D\bigl(q,Q,P_iq^{Q_i\partial_{Q_i}},
 (1-q)\partial_\alpha\bigr)\J^{\mathrm{tw}}(x,t)=0\\
 &\Longrightarrow
 D\bigl(1,Q,A_{i,\mathrm{com}}(x,t),
 \phi_\alpha\star_{x,t}\bigr)1=0.
\end{align*}
\end{thm}

\begin{proof}
Apply $D$ to
$\widetilde\J=(1-q)\widetilde S^{-1}1$ and use
\eqref{normalized-intertwining}.  The result is
$ (1-q)\widetilde S^{-1}
D(q,Q,\mathcal A_i,\phi_\alpha\star_{x,t})1=0$.
Invertibility of $\widetilde S$ and specialization at $q=1$ give
\eqref{relation-theorem-equation}.  Membership in
$\Lambda_+$ ensures that the $S$-matrix remains invertible
after specializing the inputs.  The unnormalized form is
the same calculation using
$S(P_iq^{Q_i\partial_{Q_i}})S^{-1}=\mathcal A_i$ directly.
\end{proof}

Applying the above theorem does not result in explicit relations unless the $\hat{P}_i$ can be computed. We do this using the below theorem, which we call the \emph{quantum triviality theorem}, it is a generalization of a result due to Anderson-Chen-Tseng-Iritani \cite[Lemma~6]{act}.
\begin{thm}[Quantum triviality theorem]\label{quantum-triviality}
Let $f(s_1,\dots,s_r)$ be a polynomial independent of $q$.  Suppose that,
for every $d>0$, the coefficient of $Q^d$ in
$f(P_iq^{Q_i\partial_{Q_i}})\J^{\mathrm{tw}}(x,t)/(1-q)$
vanishes at $q=\infty$.  Then
$f^{\star_{x,t}}(\widehat P_1(x,t),\dots,\widehat P_r(x,t))
=f(P_1,\dots,P_r)$.
\end{thm}

\begin{proof}
Put $B_i=P_iq^{Q_i\partial_{Q_i}}$ and
$a_f(q,Q)=f(\mathcal A_1,\dots,\mathcal A_r)1$.  From
the definitions of $\mathcal A_i$ and $\J^{\mathrm{tw}}$, we have
$f(B_1,\dots,B_r)\J^{\mathrm{tw}}/(1-q)=T_{x,t}(q)a_f(q,Q)$.
Write $T=1+\sum_{d>0}Q^dT_d$ and
$a_f=\sum_dQ^da_{f,d}$.  Proposition~2.10 of
Iritani--Milanov--Tonita \cite[Proposition~2.10]{imt} gives
$a_{f,d}\in K_T(X)\otimes\Lambda[q]$, while $T_d=O(q^{-1})$.  Induction in
effective degree gives
$[T a_f]_{Q^d}|_{q=\infty}=0\Longrightarrow a_{f,d}=0$ for $d>0$.
Consequently,
$f(A_{1,\mathrm{com}},\dots,A_{r,\mathrm{com}})1
=a_f(1,Q)=a_{f,0}=f(P_1,\dots,P_r)$.
By the commutativity of the shift operators and
\cite[Corollary~2.9]{imt}, both sides commute
with every quantum multiplication operator; hence, for every $v$,
\begin{equation*}
 f(A_{1,\mathrm{com}},\dots,A_{r,\mathrm{com}})v
 =v\star_{x,t}
 f(A_{1,\mathrm{com}},\dots,A_{r,\mathrm{com}})1
 =v\star_{x,t}f(P_1,\dots,P_r).
\end{equation*}
\end{proof}

This theorem has two important corollaries for practical computations.  Write
$\J^{\mathrm{tw}}(x,t)/(1-q)=\sum_d Q^dJ_d(q)$.

\begin{cor}\label{individual-quantum-triviality}
Fix $i$.  If $\deg_q J_d<-d_i$ for every positive degree $d$ occurring in
the preceding expansion, then $\widehat P_i(x,t)=P_i$.
\end{cor}

\begin{proof}
On the coefficient of $Q^d$, the shift
$P_iq^{Q_i\partial_{Q_i}}$ acts as multiplication by $P_iq^{d_i}$.
Thus the displayed degree bound says precisely that every
positive-degree coefficient of
$(P_iq^{Q_i\partial_{Q_i}})\J^{\mathrm{tw}}/(1-q)$ vanishes at
$q=\infty$.  Theorem~\ref{quantum-triviality}, applied to $f=s_i$, gives
$A_{i,\mathrm{com}}=P_i\star_{x,t}$.  Applying both sides to $1$ gives
$\widehat P_i=A_{i,\mathrm{com}}1=P_i\star_{x,t}1=P_i$.
\end{proof}

\begin{cor}\label{classical-quantum-products}
Assume the hypothesis of Corollary~\ref{individual-quantum-triviality}
for every $i=1,\dots,r$.  Let
$F=s_1^{c_1}\cdots s_r^{c_r}$ be a
monomial independent of $q$ such that
$\deg_q J_d<-\sum_i c_i d_i$ for all $d>0$.  Then
$F^{\star_{x,t}}(P_1,\dots,P_r)=F(P_1,\dots,P_r)$.
 
\end{cor}

\begin{proof}
Corollary~\ref{individual-quantum-triviality} gives
$F(A_{1,\mathrm{com}},\dots,A_{r,\mathrm{com}})
=F(P_1\star_{x,t},\dots,P_r\star_{x,t})$.
For a monomial $s_1^{c_1}\cdots s_r^{c_r}$, its shift operator acts on
$Q^dJ_d(q)$ by multiplication by
$P_1^{c_1}\cdots P_r^{c_r}q^{\sum_i c_i d_i}$.  The assumed bounds
therefore verify the hypothesis of Theorem~\ref{quantum-triviality}. Thus $F(A_{i,\mathrm{com}})1=F^{\star}(P_i)=F(P_i)$. 
\end{proof}

\subsection{From relations to presentations}
The above results work for obtaining relations in quantum $K$-rings, however it is also necessary to determine how many relations are sufficient to completely describe the ring, for this, we use the formal Nakayama criterion of
Gu--Mihalcea--Sharpe--Xu--Zhang--Zou
\cite[Corollary~2.7]{nakayama}. 

\begin{thm}[Formal Nakayama criterion]\label{K-Nakayama}
Let $R$ be a Noetherian integral domain, let $S$ be a Noetherian $R$-algebra,
and let $I\subset S[[Q_1,\dots,Q_r]]$ be an ideal.  Suppose that
$\varphi:S[[Q_1,\dots,Q_r]]/I\to N$
is an $R[[Q_1,\dots,Q_r]]$-algebra homomorphism, that $N$ is a finite free
$R[[Q_1,\dots,Q_r]]$-module, and that the induced map
\[
 \overline\varphi:
 \frac{S[[Q_1,\dots,Q_r]]}{I+(Q_1,\dots,Q_r)}
 \longrightarrow \frac{N}{(Q_1,\dots,Q_r)N}
\]
is an isomorphism of $R$-algebras.  Then $\varphi$ is an isomorphism.
\end{thm}

\begin{cor}\label{presentation-corollary}
Assume $K_T(X)\otimes\Lambda$ is free of finite rank.  Suppose proposed
quantum relations deform a classical presentation and their quotient is
finite free over $\Lambda$ with the same classical
basis.  If the proposed relations hold in the  twisted quantum ring,
then they give its full presentation.
\end{cor}

\subsection{Twisting Theorems}

For this section, we restrict to the symmetrized theory (i.e. put $x=0,N=1$). Given a type-$I$ twisting, Givental's quantum Adams-Riemann-Roch theorem allows us to generate a value of the big twisted $J$-function. We state this theorem in the cases we use. Suppose $E=\bigoplus_{i=1}^sE_i$ is a sum of line bundles on $X$ and write
$D_i(d)=\langle c_1(E_i),d\rangle$.
Products with an infinite lower limit are interpreted as ratios, so the
formulas are valid for positive and negative $D_i(d)$.  The following theorems are due to Givental in \cite{PermXI}. 
\begin{thm}\label{hypergeometric-twisting}
Let $\mathcal L$ denote the permutation-equivariant range of the untwisted
theory.  If $f=\sum_dQ^df_d\in\mathcal L$, then the following series belong,
respectively, to $\mathcal L^{\Eu_\lambda(E)}$,
$\mathcal L^{\Eu_\lambda(E^\vee)^{-1}}$, and
$\mathcal L^{(E,\ell)}$:
\begin{align*}
 f^{\Eu_\lambda(E)}
 &:=\sum_dQ^df_d
 \prod_{i=1}^s
 \frac{\prod_{m\leq D_i(d)}(1-\lambda E_i^{-1}q^m)}
 {\prod_{m\leq0}(1-\lambda E_i^{-1}q^m)},\\
 f^{\Eu_\lambda(E^\vee)^{-1}}
 &:=\sum_dQ^df_d
 \prod_{i=1}^s
 \frac{\prod_{m<D_i(d)}(1-\lambda E_iq^{-m})}
 {\prod_{m<0}(1-\lambda E_iq^{-m})},\\
 f^{(E,\ell)}
 &:=\sum_dQ^df_d
 \prod_{i=1}^s
 \left(E_i^{D_i(d)}q^{\binom{D_i(d)}{2}}\right)^\ell.
\end{align*}
The last line corresponds to the level-$\ell$ structure
$\det^{-\ell}(R\ft_*\ev^*E)$.
\end{thm}

\section{Reconstruction for Twisted Theories}

\label{recon}

Assume in this section that $K_T(X)$ is generated  by line bundles.

Let $\mathbf m=\{1,\dots,m\}$ and let $i\neq j$ be two named
markings.  For a partition $\mathbf m=S\sqcup S^c$ with $i\in S$ and
$j\in S^c$, and a degree splitting $d=d_1+d_2$, define
\begin{align*}
\gl_{S,d_1\mid S^c,d_2}:{}&
X_{0,S\cup\{\bullet\},d_1}\mathop{\times}_X
X_{0,S^c\cup\{\star\},d_2}
\longrightarrow X_{0,m,d},\\
D_{S,d_1\mid S^c,d_2}:={}&
\operatorname{im}(\gl_{S,d_1\mid S^c,d_2}).
\end{align*}
Only stable splittings occur in these sums.  Following
Lee--Pandharipande \cite[Theorem~1]{lp}, put
\begin{align*}
D_{i,d_1\mid j,d_2}
:={}&\sum_{\substack{S\subset\mathbf m\\i\in S,\ j\notin S}}
 D_{S,d_1\mid S^c,d_2},\\
D_{i\mid j}
 :={}&\sum_{d_1+d_2=d}D_{i,d_1\mid j,d_2}.
\end{align*}
Thus $D_{i,d_1\mid j,d_2}$ is the sum of all boundary divisors for which the
component carrying $i$ has degree $d_1$ and the component carrying $j$ has
degree $d_2$.

Let $M$ be a line bundle on $X$ and set
$\ell_e(M)=\langle c_1(M),e\rangle$.  The relations of
\cite[Theorem~1]{lp}, written additively in
$\operatorname{Pic}(X_{0,m,d})$, are
\begin{align}
 \ev_i^*M={}&\ev_j^*M+\ell_d(M)\psi_j
 -\sum_{d_1+d_2=d}\ell_{d_1}(M)
 D_{i,d_1\mid j,d_2},
 \label{LP-evaluation-relation}\\
 \psi_i+\psi_j={}&D_{i\mid j}.
 \label{LP-cotangent-relation}
\end{align}
Here $\psi_a$ denotes the additive Picard class of the cotangent line $L_a$.
Define the Cartier divisor
$B_{i\mid j}(M):=\sum_{d_1+d_2=d}\ell_{d_1}(M)D_{i,d_1\mid j,d_2}$.
Exponentiating \eqref{LP-evaluation-relation} and
\eqref{LP-cotangent-relation} give the following relations in $K(X_{g,n,d}):$
\begin{align}
 \ev_i^*M\cong{}&\ev_j^*M\otimes
 L_j^{\ell_d(M)}\otimes\cO\bigl(-B_{i\mid j}(M)\bigr),
 \label{LP-line-bundle-relation}\\
 L_i\otimes L_j\cong{}&\cO(D_{i\mid j}).
 \label{LP-line-cotangent-relation}
\end{align}

For an effective Cartier divisor $D$, recall
$1-\cO(-D)=\cO_D.$

After choosing an order on the finitely many irreducible boundary divisors
$D_\delta$ occurring in $B=\sum_\delta a_\delta D_\delta$, one has
\begin{equation}\label{multiple-divisor-telescoping}
1-\cO(-B)
=\sum_\delta
\left(\prod_{\epsilon<\delta}\cO(-a_\epsilon D_\epsilon)\right)
\left(1-\cO(-a_\delta D_\delta)\right).
\end{equation}
Every summand of
\eqref{multiple-divisor-telescoping} can be written as a line-bundle multiple of
$\cO_{D_\delta}$. 

This observation lets us prove the following reconstruction theorem:

\begin{thm}\label{mixed-reconstruction}
Let $K_T(X)$ be generated by line bundles.  Fix $N\geq1$ and a primary
permutable input $t\in\Lambda_+$.  All genus-zero
$S_n$-equivariant type-I twisted correlators with named descendant
inputs $a_i(L)$ and permutable input $t$ are determined by the corresponding
correlators with $k=1$.
\end{thm}

\begin{proof}
We induct on $k$ and the degree of a generating line bundle appearing in a coefficient of $a_i$. Given some line bundle $M$ in the $i$th factor, we can use \ref{LP-line-bundle-relation} to move the input to another factor, lowering the degree in $M$, the cost is adding a descendent term and the term $\mathcal{O}(-B_{i|j}(M))$, however, this can be rewritten as $1-\mathcal{O}_{B_{i|j}(M)}$ and split into two terms by linearity. The first term has a smaller degree in $M$ in the $i$th factor, and the second term can be computed via the splitting axiom, which expresses it as a sum of invariants with fewer unpermuted marked points. 

If we have two marked points with nontrivial descendent inputs, we can use \ref{LP-line-cotangent-relation} to remove one descendent degree at each input, and apply the splitting axiom on the resulting boundary divisor. 

A primary input that has no power of the generating line bundle is just an element of $\Lambda$, so we can apply the string equation, which reduces $k$.

\end{proof}

\section{Quantum Abelian-Non/Abelian Correspondence}
\label{nonab}
\subsection{Classical correspondence}

Throughout this section, let $G$ be a connected complex reductive algebraic
group, let $T\subset G$ be a maximal algebraic torus, and let
$W=N_G(T)/T$ be the Weyl group.  Let $V$ be a smooth quasiprojective complex
$G$-variety with a linearization for which the GIT quotients
$Z:=V\sslash G$ and $Y:=V\sslash T$
are smooth and projective.  We assume stable and semistable loci agree and
that the unstable locus has codimension at least two. 

Let $\mathcal{R}$ denote the bundle on $Y$ given by the sum of the line
bundles associated to the roots, and define
$\operatorname{ann}_{K(Y)^W}(\Eu(\mathcal R))
:=\{a\in K(Y)^W:a\Eu(\mathcal R)=0\text{ in }K(Y)\}$.  The specialization
of Harada--Landweber's more general theorem used here is
\begin{equation*}
 K(Z)\cong
 \frac{K(Y)^W}
 {\operatorname{ann}_{K(Y)^W}(\Eu(\mathcal R))}.
\end{equation*}
We denote the quotient map by
$\phi:K(Y)^W\twoheadrightarrow K(Z)$.

We also have, by \cite[Theorem C]{Martin}:
$\chi_Z(\phi(a))=|W|^{-1}\chi_Y(a\Eu(\mathcal R))$.

\subsection{Quantum formulation}

Following the convention introduced earlier, let $\mathcal L_Z$ denotes the
permutation-equivariant range of $t\mapsto\J_Z(0,t)$, and let
$\mathcal L_Y^{tw}$ denote the range of the
$\Eu_\lambda(\mathcal R)$-twisted theory of $Y$.  For fixed permutable
inputs $t$ on $Z$ and $\widetilde t$ on $Y$, the Lagrangian
cones at permutable input $t$ are denoted by $\mathcal L_{Z,t}$ and
$\mathcal L^{tw}_{Y,\widetilde t}$, respectively.

For a vector bundle $M$ on $Z$ obtained from a $G$-representation let $\widetilde M$ be the vector bundle bundle on $Y$ obtained by restricting $G$ to $T$.  Compatibility of the Kirwan maps gives
$\phi(\widetilde M)=M$.  Define the specialization of effective curve classes
$\rho$ to be dual to this; explicitly,
\begin{align*}
 \rho&:H_2(Y)\longrightarrow H_2(Z),\\
 \left\langle c_1(M),\rho(\widetilde d)\right\rangle_Z
 &=\left\langle c_1(\widetilde M),\widetilde d\right\rangle_Y
 \quad\text{for every such }M\text{ and }\widetilde d\in H_2(Y).
\end{align*}
This induces a specialization of Novikov variables.  If
$F=\sum_{\widetilde d}Q^{\widetilde d}F_{\widetilde d}$ has $W$-invariant
orbit sums, define
\begin{equation*}
 \Phi_Q(F):=
 \sum_dQ^d\lim_{\lambda\to1}\phi\left(
 \sum_{\rho(\widetilde d)=d}F_{\widetilde d}\right).
\end{equation*}

\begin{conj}
\label{-ana-cone}\label{germ}\label{sgerm}
For every pair of corresponding permutable inputs $\widetilde t,t$,
coefficientwise application of $\Phi_Q$ identifies a $W$-invariant formal
neighborhood of $1-q$ in
$\mathcal L^{tw}_{Y,\widetilde t}$ with a formal neighborhood of
$1-q$ in $\mathcal L_{Z,t}$.  It respects positive-space projections and
maps each source ruling to the corresponding target ruling; every target
ruling has a local $W$-invariant lift. 
\end{conj}

\begin{conj}
\label{-ana-correlators}\label{cor}\label{corsym}
Let $a_1(L),\dots,a_k(L)$ and $t(L)$ be descendant inputs on $Z$, and choose
$W$-invariant lifts $\widetilde a_i(L),\widetilde t(L)$ on $Y$.  For every
$g,k,n,d$:
\begin{align}
 &\left\langle a_1(L),\dots,a_k(L);t(L),\dots,t(L)
 \right\rangle^{S_n}_{g,k+n,d;Z}\notag\\
 &\quad=\frac1{|W|}\lim_{\lambda\to1}
 \sum_{\rho(\widetilde d)=d}
 \left\langle
 \widetilde a_1(L),\dots,\widetilde a_k(L);
 \widetilde t(L),\dots,\widetilde t(L)
 \right\rangle^{S_n,tw,\lambda}_{g,k+n,\widetilde d;Y}.
 \label{-correspondence-formula}
\end{align}
In particular, the right-hand side has a well-defined limit and is
independent of the chosen lifts.
\end{conj}

\begin{conj}
\label{-ana-rings}\label{ring}\label{ringsym}
For corresponding $W$-invariant primary inputs
$(\widetilde x,\widetilde t)$ on $Y$ and $(x,t)$ on $Z$, and for every
$N\geq1$, the map $\Phi_Q$ is a surjective homomorphism
\begin{equation*}
 \Phi_Q:
 \bigl(QK^{tw,S_n}(Y),
 \star_{\widetilde x,\widetilde t}\bigr)^W
 \longrightarrow
 \bigl(QK^{S_n}(Z),\star_{x,t}\bigr).
\end{equation*}

The source is specialized at $\lambda=1$ only after taking the twisted
products and the degree sums in the definition of $\Phi_Q$.  At
$\widetilde t=t=0$ this is the ordinary ring statement; at
$\widetilde x=x=0$, followed by extraction of the trivial
$S_n$-representation coefficient, it is the symmetrized ring statement.
\end{conj}

\subsection{Relations among the three formulations}

\begin{prop}\label{j}
Conjecture~\ref{-ana-cone} implies $\Phi_Q(J_Y^{tw})=J_Z$.
More generally it identifies the $J$-functions at corresponding
inputs.
\end{prop}

\begin{proof}
Take $x=t=0$. 
\end{proof}

\begin{thm}\label{correlators-imply-rings}
Conjecture~\ref{-ana-correlators} implies
Conjecture~\ref{-ana-rings}.
\end{thm}

\begin{proof}

With two named primary insertions, \eqref{-correspondence-formula}
says that $\rho_*G_Y^{tw}(a,b)$ tends to
$|W|G_Z(\Phi_Qa,\Phi_Qb)$ as $\lambda\to1$.
The same formula with three insertions and the definition of the products
give, for every lifted test class $c$,
\begin{align*}
 G_Z\bigl(\Phi_Q(a\star_Yb),\Phi_Qc\bigr)
 &=\frac1{|W|}\lim_{\lambda\to1}
 \rho_*\llangle a,b,c\rrangle_Y^{tw}\\
 &=\llangle\Phi_Qa,\Phi_Qb,\Phi_Qc\rrangle_Z\\
 &=G_Z\bigl(\Phi_Qa\star_Z\Phi_Qb,\Phi_Qc\bigr).
\end{align*} 
\end{proof}

\begin{prop}\label{cone-implies-one-point}
Conjecture~\ref{-ana-cone}, together with the classical integration formula
above, implies the genus-zero one-named-point case of
Conjecture~\ref{-ana-correlators} for every permutable input.
\end{prop}

\begin{proof}
Fix $N\geq1$.  For $c\in K(Z)$ and a $W$-invariant lift
$\widetilde c\in K(Y)$, put
\begin{align*}
 \mathcal C^Z_{d,n}(c;x,t;q)
 &:=\sum_{k\geq0}\frac1{k!}
 \left\langle\frac{c}{1-qL},x,\dots,x;
 t,\dots,t\right\rangle^{S_n}_{0,1+k+n,d;Z},\\
 \mathcal C^Y_{\widetilde d,n}(\widetilde c;\widetilde x,
 \widetilde t;q)
 &:=\sum_{k\geq0}\frac1{k!}
 \left\langle\frac{\widetilde c}{1-qL},
 \widetilde x,\dots,\widetilde x;
 \widetilde t,\dots,\widetilde t
 \right\rangle^{S_n,tw,\lambda}_{0,1+k+n,\widetilde d;Y}.
\end{align*}

The cone correspondence gives
$\Phi_Q(\J_Y^{tw}(\widetilde x,\widetilde t))=\J_Z(x,t)$.
Pairing this equality with $c$ and applying the classical integration
formula gives
\begin{align*}
 &\bigl(\J_Z(x,t)-(1-q)-x-t,c\bigr)_Z\\
 &=\frac1{|W|}\lim_{\lambda\to1}\rho_*
 \bigl(\J_Y^{tw}(\widetilde x,\widetilde t)
 -(1-q)-\widetilde x-\widetilde t,\widetilde c
 \bigr)^{tw}_Y.
\end{align*}
Taking the coefficient of $Q^d$ and the homogeneous contribution with $n$
copies of $U_N$ yields
$\mathcal C^Z_{d,n}(c;x,t;q)=|W|^{-1}\lim_{\lambda\to1}
\sum_{\rho(\widetilde d)=d}\mathcal C^Y_{\widetilde d,n}
(\widetilde c;\widetilde x,\widetilde t;q)$,
which is \eqref{-correspondence-formula} for one named
descendant insertion.
\end{proof}

\begin{thm}\label{imp2}\label{imp2sym}
Assume the $K$-rings of $Y$ and $Z$ are
generated by line bundles.  Then Conjecture~\ref{-ana-cone} implies the
genus-zero correlator correspondence.  More precisely, fix $k,n,d$ and
$N\geq1$, named descendant inputs $a_i(L)$ and a primary permutable
input $t$ on $Z$, and $W$-invariant lifts $\widetilde a_i(L)$ and
$\widetilde t$ on $Y$.
 \begin{align*}
 &\left\langle a_1(L),\dots,a_k(L);t,\dots,t
 \right\rangle^{S_n}_{0,k+n,d;Z}\\
 &\quad=\frac1{|W|}\lim_{\lambda\to1}
 \sum_{\rho(\widetilde d)=d}
 \left\langle
 \widetilde a_1(L),\dots,\widetilde a_k(L);
 \widetilde t,\dots,\widetilde t
 \right\rangle^{S_n,tw,\lambda}_{0,k+n,\widetilde d;Y}.
\end{align*}
\end{thm}

\begin{proof}
Proposition~\ref{cone-implies-one-point} proves the assertion for one
named insertion.  Starting from this case, run the reconstruction
induction of Theorem~\ref{mixed-reconstruction} simultaneously for $Y$ and
$Z$.  If $\widetilde M$ and $M$ are corresponding line bundles and
$\rho(\widetilde d)=d$, then
$\phi(\widetilde M)=M$ and
$\langle c_1(\widetilde M),\widetilde d\rangle
=\langle c_1(M),d\rangle$.
Consequently the line-bundle and cotangent-line relations
have the same coefficients on the two sides.  Linear combinations, the
Novikov specialization, and the string equation also preserve the
correspondence.  It remains to prove that every use of the splitting axiom in
the reconstruction preserves it.

\smallskip
\noindent\emph{Compatibility of splitting.}
Fix one separating boundary term, with degree decomposition $d=d_1+d_2$ and
with $n_1$ and $n_2$ permutable markings on its two components.  Leave the
insertion at the node open.  Denote the resulting component covectors on
$K(Z)$ by $A_Z$ and $B_Z$, and denote by $A_Y^\lambda$ and $B_Y^\lambda$ the
corresponding $\Eu_\lambda(\mathcal R)$-twisted covectors after summing the
first component over $\rho(\widetilde d_1)=d_1$ and the second over
$\rho(\widetilde d_2)=d_2$.  Suppose, as in the reconstruction induction,
that the correspondence is already known for these component correlators.
Thus, for every $W$-invariant class $\widetilde c$,
\begin{align*}
 \lim_{\lambda\to1}A_Y^\lambda(\widetilde c)
 &=|W|A_Z(\phi(\widetilde c)),&
 \lim_{\lambda\to1}B_Y^\lambda(\widetilde c)
 &=|W|B_Z(\phi(\widetilde c)).
\end{align*}
We claim that the correspondence then holds for the boundary term obtained by
contracting these covectors at the node.

Let
$g_Z(a,b):=\chi_Z(ab)$ and
$g_Y^{tw}(\widetilde a,\widetilde b)
:=\chi_Y(\widetilde a\widetilde b\Eu_\lambda(\mathcal R))$.
The classical integration formula gives, at $\lambda=1$,
\begin{equation*}
 g_Y^{tw,1}(\widetilde a,\widetilde b)
 =|W|g_Z(\phi(\widetilde a),\phi(\widetilde b)).
\end{equation*}
Since $\ker\phi$ is the annihilator of
$\Eu(\mathcal R)$ in $K(Y)^W$,  the twisted pairing therefore descends to the
nondegenerate quotient $K(Y)^W/\ker\phi$, identified with $K(Z)$.  The
degree-summed source covectors are
$W$-invariant, and their limits factor through $\phi$ by the induction
hypothesis.  For a basis $\{\phi_\mu\}$ of $K(Z)$ and
$W$-invariant lifts $\{\widetilde\phi_\mu\}$, its inverse matrix therefore
satisfies
$(g_Y^{tw,1})^{\mu\nu}=|W|^{-1}(g_Z)^{\mu\nu}$ on this quotient.
The factor appearing at the node on the space of maps to $Y$ in \eqref{twist-separating} is
$\Eu_\lambda(\mathcal R)^{-1}$, so the splitting axiom contracts the $Y$
covectors with $(g_Y^{tw})^{-1}$; the contractipon on $Z$ uses
$g_Z^{-1}$.  The induction hypothesis now gives:
\begin{align*}
 &\frac1{|W|}\lim_{\lambda\to1}
 \sum_{\mu,\nu}A_Y^\lambda(\widetilde\phi_\mu)
 (g_Y^{tw})^{\mu\nu}
 B_Y^\lambda(\widetilde\phi_\nu)\notag\\
 &\quad=\frac1{|W|}\sum_{\mu,\nu}
 (|W|A_Z(\phi_\mu))
 \bigl(|W|^{-1}(g_Z)^{\mu\nu}\bigr)
 (|W|B_Z(\phi_\nu))\notag\\
 &\quad=\sum_{\mu,\nu}A_Z(\phi_\mu)(g_Z)^{\mu\nu}B_Z(\phi_\nu),
\end{align*}
This is exactly the contribution from the boundary divisor on $Y$. 
It remains only to check that this contraction is the splitting formula for
the $S_n$-equivariant invariants.  If $V_1$ and $V_2$ are the virtual
$S_{n_1}$- and $S_{n_2}$-modules of the two components, then the module on the
boundary is
$\operatorname{Ind}_{S_{n_1}\times S_{n_2}}^{S_n}(V_1\boxtimes V_2)$.
After inserting $U_N$ at the permutable markings,
\eqref{equivariant-binomial-expanded} gives
\begin{align*}
 &\left[
 \operatorname{Ind}_{S_{n_1}\times S_{n_2}}^{S_n}(V_1\boxtimes V_2)
 \otimes U_N^{\otimes n}
 \right]^{S_n}\\
 &=\left[V_1\otimes U_N^{\otimes n_1}\right]^{S_{n_1}}
 \otimes
 \left[V_2\otimes U_N^{\otimes n_2}\right]^{S_{n_2}}.
\end{align*}
Thus the two factors in the contraction are precisely the two component
correlators.  The induced $S_n$-module already accounts for every choice of
the $n_1$ permutable markings on the first component, so there is no
additional binomial coefficient.  Additivity of $\rho$ gives
\begin{align*}
 \sum_{\substack{\rho(\widetilde d)=d\\
 \widetilde d_1+\widetilde d_2=\widetilde d}}(\cdots)
 =\sum_{d_1+d_2=d}
 \sum_{\substack{\rho(\widetilde d_1)=d_1\\
 \rho(\widetilde d_2)=d_2}}(\cdots),
\end{align*}
so the degree sums in the two splitting formulas also agree.  This proves the
compatibility statement for every separating boundary term.

\end{proof}

\begin{rem}
For Grassmannians and flag varieties the nonequivariant $K$-ring need not be
generated by line bundles in the required form.  One may first work
equivariantly and localize near the augmentation ideal, where fixed-point
localization supplies line-bundle generators, perform reconstruction, and
then take the nonequivariant limit when it exists.
\end{rem}

\subsection{Additional twistings}

Let $E$ be a $G$-representation inducing bundles $E_Z$ and $E_Y$ on the two
quotients, and let $C$ be an invertible multiplicative class.  Each of the
three conjectures has a simultaneously twisted version: twist the target by
$C(E_Z)$ and the source by $\Eu_\lambda(\mathcal R)C(E_Y)$.
The statements and implications above are unchanged because restriction,
splitting, and the permutation-equivariant invariants commute with both
twisting classes.

Givental's Quantum Adams--Riemann--Roch theorem
\cite[Theorem~1]{PermXI} expresses the additional twistings by loop-group
operators $\Delta_Z$ and $\Delta_Y$.  Compatibility of the associated
bundles gives $\Phi_Q\circ\Delta_Y=\Delta_Z\circ\Phi_Q$.
Consequently the untwisted cone correspondence implies every such
compatibly twisted cone correspondence, provided the indicated
nonequivariant limits exist.  

\subsection{Computing rings by abelianization}

Assume Conjecture~\ref{-ana-rings}.  Let
$P_1,\dots,P_r$ be line-bundle generators of $K(Y)$ dual to Novikov
coordinates, and choose a $W$-invariant
$I\in\mathcal L^{tw}_{Y,\widetilde t}$ satisfying
$\Phi_Q(I)=J_Z$.  The series $I$ need not be the twisted small $J$-function
of $Y$.  By
\cite[Theorem~2 and Corollary~1]{PermVII}, it can be
written uniquely near $1-q$ as $I=(1-q)T_\tau(q)h(q)$ for a primary input
$\tau$, where $h\in\Kloop_+$ and
$h(1)-1\in\Lambda_+$.

We further have:
\begin{thm}
$h$ and $\tau$ are $W$-invariant, $\Phi_Q(\tau)=0$, and $\Phi_Q(h)=1$.
\end{thm}
\begin{proof}
The $W$-equivariant formal implicit-function theorem makes both $\tau$ and
$h$ Weyl invariant.  Since $\Phi_Q$ respects rulings,
$\Phi_Q(\tau)=0$ and $\Phi_Q(h)=1$.
\end{proof}

\begin{thm}
\label{usable}\label{abelianized-relations}
Let $D$ be a polynomial finite-difference operator annihilating the
logarithmically normalized form of $I$.  Then
$D^{\star_Y}(1,Q,\widehat P_i(\tau))=0$ in $QK^{tw,S_n}(Y)$ at $\tau$.
\end{thm}

\begin{proof}
Intertwining the shifts with $T_\tau$ rewrites $DI=0$ as
$(1-q)T_\tau D(q,Q,\mathcal A_i)h(q)=0$.  At $q=1$ this becomes
$D(1,Q,A_{i,\mathrm{com}}(\tau))h(1)=0$.
By \eqref{commuting-shift-multiplication}, this endomorphism is quantum
multiplication by the class
$D^{\star_Y}(1,Q,\widehat P_i(\tau))$.  Since
$h(1)-1\in\Lambda_+$, the class $h(1)$ is a unit for
$\star_Y$, and therefore
$0=D^{\star_Y}(1,Q,\widehat P_i(\tau))\star_Yh(1)$ implies
$D^{\star_Y}(1,Q,\widehat P_i(\tau))=0$.
\end{proof}

There is also a form of quantum triviality adapted to the correspondence.  Put
$\mathbf K_Y:=K(Y)\otimes\Lambda$ and
$\mathbf K_Z:=K(Z)\otimes\Lambda$, and write
$\star_Y:=\star^{tw}_{Y,\tau}$ for the product in
$QK^{tw,S_n}(Y)$ at the input $\tau$. 
Regard
$\Phi_Q:\mathbf K_Y^W\to\mathbf K_Z$ as the $\Lambda$-linear map of
modules defined coefficientwise above, and set
$K_\Phi:=\ker_\Lambda(\Phi_Q:\mathbf K_Y^W\to\mathbf K_Z)$. 
The definition of $\Phi_Q$ also respects the ordinary tensor-product
multiplication on $W$-invariant coefficients, so $K_\Phi$ is a classical
ideal in $\mathbf K_Y^W$.  We define
\begin{align*}
 \mathscr I
 &:=\left\langle K_\Phi\right\rangle_{\mathbf K_Y,\mathrm{cl}}
   =\mathbf K_Y\cdot_{\mathrm{cl}}K_\Phi\\
 &:=\left\{
 \sum_{\nu=1}^N a_\nu\cdot_{\mathrm{cl}} b_\nu:
 a_\nu\in\mathbf K_Y,\quad b_\nu\in K_\Phi
 \right\}.
\end{align*}
Thus $\mathscr I$ is generated by the \emph{module kernel} of $\Phi_Q$ as
an ideal in the entire ring $\mathbf K_Y$, using classical multiplication.

Since $|W|$ is invertible,
$\mathscr I\cap\mathbf K_Y^W=K_\Phi$.
Indeed, if $u=\sum_\nu a_\nu\cdot_{\mathrm{cl}}b_\nu$ is $W$-invariant,
then, because every $b_\nu$ is invariant,
\begin{equation*}
 u=\frac1{|W|}\sum_{w\in W}w(u)
  =\sum_\nu\left(\frac1{|W|}\sum_{w\in W}w(a_\nu)\right)
    \cdot_{\mathrm{cl}}b_\nu\in K_\Phi.
\end{equation*}

Write $I/(1-q)=\sum_{\widetilde d}Q^{\widetilde d}I_{\widetilde d}(q)$.
For $u,v\in\mathbf K_Y$, the notation
$u\equiv v\pmod{\mathscr I}$ means $u-v\in\mathscr I$ in the underlying
$\Lambda$-module.  For endomorphisms $A,B$ of $\mathbf K_Y$, it means the
explicit coefficientwise assertion
$(A-B)(\mathbf K_Y)\subset\mathscr I$ in that underlying module.  The
commuting specializations $A_{i,\mathrm{com}}(\tau)$
are endomorphisms, whereas $\widehat P_i(\tau)$ are classes, and
\eqref{commuting-shift-multiplication} reads
$A_{i,\mathrm{com}}(\tau)v=v\star_Y\widehat P_i(\tau)$ for
$v\in\mathbf K_Y$.
Moreover, the $W$-invariance of $h$ and the  equality
$\Phi_Q(h)=1$ imply that every coefficient of $h(q)-1$ belongs to
$K_\Phi\subset\mathscr I$.  Thus, 
$h(q)\equiv1\pmod{\mathscr I}$ and
$h(1)\equiv1\pmod{\mathscr I}$.

\begin{thm}
\label{abelianized-triviality}
Let $F(s_1,\dots,s_r)$ be any polynomial independent of $q$.  If every
positive-degree coefficient of
$F(P_iq^{Q_i\partial_{Q_i}})I/(1-q)$
vanishes at $q=\infty$ modulo $\mathscr I$ in the underlying coefficient
module $\mathbf K_Y$, then
$F^{\star_Y}(\widehat P_1(\tau),\dots,\widehat P_r(\tau))
\equiv F(P_1,\dots,P_r)\pmod{\mathscr I}$.
\end{thm}

\begin{proof}
Put $B_i=P_iq^{Q_i\partial_{Q_i}}$.  Intertwining the shifts with
 $T_\tau$ in the expression $I=(1-q)T_\tau h$ gives
\begin{equation*}
 F(B_1,\dots,B_r)\frac{I}{1-q}
 =T_\tau(q)a_F(q,Q),
 a_F(q,Q):=F(\mathcal A_1,\dots,\mathcal A_r)h(q).
\end{equation*}
Since $h(q)\equiv1\pmod{\mathscr I}$,
$a_F(q,Q)\equiv F(\mathcal A_1,\dots,\mathcal A_r)1\pmod{\mathscr I}$.
Its coefficients are polynomial in $q$ by
\cite[Proposition~2.10]{imt}, whereas every positive-degree coefficient of
$T_\tau(q)-1$ is $O(q^{-1})$.  Induction in effective degree in
the preceding identity therefore gives
$[a_F]_{Q^{\widetilde d}}=0$ modulo $\mathscr I$ in $\mathbf K_Y[q]$ for
$\widetilde d>0$.
The degree-zero terms are $h_0=1$ and $\mathcal A_{i,0}=P_i$.  Specializing
at $q=1$ consequently yields
\begin{equation*}
 F(A_{1,\mathrm{com}}(\tau),\dots,A_{r,\mathrm{com}}(\tau))h(1)
 \equiv F(P_1,\dots,P_r)\pmod{\mathscr I}.
\end{equation*}
Since $h(1)\equiv1\pmod{\mathscr I}$ and
$F(A_{1,\mathrm{com}}(\tau),\dots,A_{r,\mathrm{com}}(\tau))$
preserves $\mathscr I$, this gives
\begin{equation*}
 F(A_{1,\mathrm{com}}(\tau),\dots,A_{r,\mathrm{com}}(\tau))1
 \equiv F(P_1,\dots,P_r)\pmod{\mathscr I}.
\end{equation*}
Finally, for every $v\in\mathbf K_Y$,
\begin{align*}
 &\bigl(F(A_{1,\mathrm{com}}(\tau),\dots,A_{r,\mathrm{com}}(\tau))
       -F(P_1,\dots,P_r)\star_Y\bigr)v\\
 &=v\star_Y\bigl(
 F(A_{1,\mathrm{com}}(\tau),\dots,A_{r,\mathrm{com}}(\tau))1
 -F(P_1,\dots,P_r)\bigr)
 \in\mathscr I,
\end{align*}
\end{proof}

\begin{cor}\label{coras}
Fix $i$.  If
 $\deg_q I_{\widetilde d}<-\widetilde d_i$
for every positive degree contributing after specialization, then
$\widehat P_i(\tau)\equiv P_i\pmod{\mathscr I}$
as classes in the underlying module.
\end{cor}

\begin{proof}
On the coefficient of $Q^{\widetilde d}$, the shift
$P_iq^{Q_i\partial_{Q_i}}$ multiplies by $P_iq^{\widetilde d_i}$.  The
displayed bound is therefore exactly the hypothesis of
Theorem~\ref{abelianized-triviality} for $F=s_i$.
\end{proof}

\begin{cor}\label{corbs}
Assume the hypothesis of Corollary~\ref{coras} for every $i=1,\dots,r$.
Let $F=\sum_{\mathbf c}f_{\mathbf c}s_1^{c_1}\cdots s_r^{c_r}$ be a
polynomial such that, for every monomial with $f_{\mathbf c}\neq0$,
$\deg_q I_{\widetilde d}<-\sum_i c_i\widetilde d_i$
for every positive degree contributing after specialization.  Then:
$F^{\star_Y}(P_1,\dots,P_r)
\equiv F(P_1,\dots,P_r)\pmod{\mathscr I}$.
If $F$ is $W$-invariant, we further have:
$\Phi_QF^{\star}
(P_1,\dots,P_r)=\Phi_QF(P_1,\dots,P_r)$.
\end{cor}

\begin{proof}
Corollary~\ref{coras} and the fact that all operators involved preserve
$\mathscr I$ give
$F(A_{1,\mathrm{com}}(\tau),\dots,A_{r,\mathrm{com}}(\tau))
\equiv F(P_1\star_Y,\dots,P_r\star_Y)\pmod{\mathscr I}$.
For a monomial $s_1^{c_1}\cdots s_r^{c_r}$, applying the corresponding
shift to $Q^{\widetilde d}I_{\widetilde d}(q)$ raises its $q$-degree by
$\sum_i c_i\widetilde d_i$.  The assumed bounds therefore let us apply
Theorem~\ref{abelianized-triviality} monomial by monomial, giving
$F(A_{1,\mathrm{com}}(\tau),\dots,A_{r,\mathrm{com}}(\tau))
\equiv F(P_1,\dots,P_r)\star_Y\pmod{\mathscr I}$.
Comparing these congruences and applying both sides to $1$ proves
the first assertion.  If $F$ is $W$-invariant,
the difference belongs to
$\mathscr I\cap\mathbf K_Y^W=K_\Phi$, which proves the second assertion.
\end{proof}

\section{Applications}
\label{apps}
\subsection{Level structures on projective space}

We first consider $X=\mathbb{P}^N$ with the level structure twisting determined by the bundle $P:=O(-1)$ and the integer $\ell$. The small $J$-function of this theory, for $-1<\ell\leq N+1$, was calculated by Givental-Yan in \cite{Xiaohan} and is equal to:

$$J^\ell:=(1-q)\sum_d Q^d\frac{P^\ell q^{\binom{d}{2}}}{\prod_i\prod_{j=1}^d(1-Pq^j/\Lambda_i)}$$

\begin{rem}
    The terms $\frac{1}{(1-Pq^j/\Lambda_i)}$ is handled in the following way. It is rewritten using the expansion

    \begin{equation}
        \label{expansion}
        \frac{1}{(1-Pq^j/\Lambda_i)}=\sum_{n} \frac{q^{jn}(P/\Lambda_i-1)^n}{(1-q^j)^n}=\sum_{n} \frac{q^{jn}((P-1)(\Lambda_i^{-1}-1)+(P-1)+(\Lambda_i^{-1}-1)-1)^n}{(1-q^j)^n}
    \end{equation}

    The right hand expression is well-defined after completing with respect to $\Lambda^{-1}-1$, and is a rational function with poles at roots of unity with $q$-degree $-j$. 
\end{rem}
Now, we can find 
\begin{prop}
  $\tilde{J}^\ell$ satisfies the following difference equation:
  $$\prod_i(1-\Lambda_i^{-1}q^{Q\p_{Q}})\tilde{J}^\ell=Qq^{\ell Q\p_{Q}} J^\ell_{1,N+1}$$

\end{prop}

\begin{proof}
\begin{equation}
\begin{split}
&\prod_i(1-\Lambda_i^{-1}q^{Q\p_{Q}})P^{\frac{\ln(Q)}{\ln(q)}}J^{\ell}_{1,N+1}=(1-q)\prod_i(1-q^dP/\Lambda_i)P^{\frac{\ln(Q)}{\ln(q)}}\sum_{d\geq 0}\frac{Q^d P^{\ell d}q^{\ell\binom{d}{2}}}{\prod_i\prod_{m=1}^{d}(1-q^mP/\Lambda_i)}\\
=&(1-q)P^{\frac{\ln(Q)}{\ln(q)}}\sum_{d\geq 1}\frac{Q^dP^{\ell d}q^{\ell\binom{d}{2}}}{\prod_i\prod_{m=1}^{d-1}(1-q^mP/\Lambda_i)}=
(1-q)P^{\frac{\ln(Q)}{\ln(q)}}\sum_{d\geq 1}q^{\ell (d-1)}P^{\ell d}\frac{Q^dP^{\ell (d-1)}q^{\ell\binom{d-1}{2}}}{\prod_i\prod_{m=1}^{d-1}(1-q^mP/\Lambda_i)}\\
=&(1-q)Qq^{\ell Q\p_{Q}}P^{\frac{\ln(Q)}{\ln(q)}}\sum_{d\geq 1}\frac{Q^{d-1}P^{\ell (d-1)}q^{\ell\binom{d-1}{2}}}{\prod_i\prod_{m=1}^{d-1}(1-q^mP/\Lambda_i)}=Qq^{\ell Q\p_{Q}}\tilde{J}_{1,N+1}
\end{split}
\end{equation}

\end{proof}

Writing $J^\ell$ as $(1-q)\sum Q^dJ_d$, we have that $deg(J_d)=\ell \binom{d}{2}-(N+1)\binom{d+1}{2}$, this is always strictly less than $-d$, so by the quantum triviality theorem $\hat{P}=P$. So we can conclude:

\begin{cor}
    The quantum $K$-ring with level structure of $\mathbb{P}^N$ with respect to the bundle $P$ and the level $\ell$ s determined by the relation:

    \begin{equation}\prod_i (1-P/\Lambda_i)=QP^\ell\end{equation}
\end{cor}

\subsection{Complete intersections in projective space}

Let $i:Y\hookrightarrow\mathbb P^N$ be a smooth complete intersection of
multidegree $(\ell_1,\dots,\ell_r)$, and continue to write
$P=\cO(-1)$. Further assume $\ell_i>1$, otherwise we can express the complete intersection in a smaller projective space. Put
$E=\bigoplus_{i=1}^r\cO(\ell_i)=\bigoplus_{i=1}^rP^{-\ell_i}$.
Quantum Lefschetz \cite[Theorem~3]{PermXI} states that the following function is in the range of the $S_n$-equivariant $Eu_\lambda(E)$-twisted $J$-function: 
\begin{equation*}
 J^{E,\lambda}=(1-q)\sum_{d\geq0}Q^d
 \frac{\displaystyle
 \prod_{i=1}^r\prod_{m=1}^{\ell_i d}
 (1-\lambda P^{\ell_i}q^m)}
 {\displaystyle\prod_{m=1}^{d}(1-Pq^m)^{N+1}}.
\end{equation*}

If $J^{E,\lambda}=(1-q)\sum_dQ^dJ_d$, then define
\begin{equation*}
 \delta_d:=\deg_q J_d
 =\sum_{i=1}^r\binom{\ell_i d+1}{2}
 -(N+1)\binom{d+1}{2}
 =\frac d2\left[
 d\Bigl(\sum_i\ell_i^2-(N+1)\Bigr)
 +\sum_i\ell_i-(N+1)\right].
\end{equation*}

\begin{prop}\label{complete-intersection-bound}
Suppose $\sum_i\ell_i^2\leq N+1$ and
$\sum_i(\ell_i^2+\ell_i)\leq2N-1$.
Then $\delta_d<-d$ for every $d>0$.  Hence
$J^{E,\lambda}$ is the small Euler-twisted
$J$-function and $\widehat P=P$. 

If $N>3$, the second condition is automatically satisfied. 

\end{prop}

\begin{proof}
Set $A=\sum_i\ell_i^2-(N+1)$ and
$B=\sum_i\ell_i-(N+1)$.  By
the formula for $\delta_d$, the inequality $\delta_d<-d$ is equivalent to
$Ad+B+2<0$.  The first inequality in
the hypothesis gives $A\leq0$, so it is
enough to check $d=1$; the second inequality is precisely $A+B+2<0$.
This decay also makes the positive-space projection equal to $1-q$.

For the final assertion, $\ell_i\geq2$ gives
$\sum_i\ell_i\leq\frac12\sum_i\ell_i^2\leq(N+1)/2$.
When $N>3$, this implies the second inequality.  
\end{proof}

The coefficient recurrence for $J^{E,\lambda}$ is
equivalent to
\begin{equation}\label{complete-intersection-difference}
 (1-Pq^{Q\partial_Q})^{N+1}J^{E,\lambda}
 =Q\prod_{i=1}^r\prod_{m=1}^{\ell_i}
 \left(1-\lambda q^m(Pq^{Q\partial_Q})^{\ell_i}\right)
 J^{E,\lambda}.
\end{equation}

Let $QK^{amb}(Y)$ be the ambient part of $QK(Y)$, i.e. both elements of $K(Y)$ and Novikov parameters are those pulled back from $\mathbb{P}^N$.
\begin{thm}\label{complete-intersection-relation}
Under the hypotheses of Proposition~\ref{complete-intersection-bound}, the
following relation holds in the Euler-twisted quantum ring of
$\mathbb P^N$:
$(1-P)^{N+1}=Q\prod_{i=1}^r(1-\lambda P^{\ell_i})^{\ell_i}$.
After taking $\lambda\to1$ and restricting to $Y$, one obtains
$(1-P|_Y)^{N+1}=Q\prod_{i=1}^r(1-P|_Y^{\ell_i})^{\ell_i}$ in $QK^{amb}(Y)$.
\end{thm}

\begin{proof}
By quantum triviality, we can replace $\widehat P$ by $P$.  Quantum Lefschetz identifies
the nonequivariant restriction of the twisted theory with the ambient part
of $QK(Y)$.  
\end{proof}

For a hypersurface this proves, in the stated numerical range, the
$K$-theoretic relation predicted by the physical analysis of
\cite[\S5.1]{gu2022symp}.
\subsection{Quasimaps to $T^*\mathbb P^N$}

Pushkar--Smirnov--Zeitlin define a quasimap quantum $K$-ring for cotangent
bundles of Grassmannians and quantum deformations of their tautological
bundles \cite[Theorem~1]{Bax}.  For
$T^*\mathbb P^N=T^*\operatorname{Gr}(1,N+1)$, they obtain:
$\prod_{j=0}^{N}(s-\Lambda_j)/(\hbar\Lambda_j-s)
=z\hbar^{-(N+1)/2}$.
The class $s$ is the eigenvalue of what the authors refer to as the \emph{quantum tautological line bundle}, a $z$-deformation of $\mathcal{O}(1)$

On the stable-map side, introduce $\lambda$ for the equivariant parameter scaling the fibers of $T^*\mathbb{P}^N$ and consider
the series
\begin{equation*}
 \widehat I_N=(1-q)\sum_{d\geq0}Q^d
 \prod_{j=0}^{N}
 \frac{\prod_{m=0}^{d-1}(1-\lambda P\Lambda_j^{-1}q^m)}
 {\prod_{m=1}^{d}(1-P\Lambda_j^{-1}q^m)}.
\end{equation*}
By 
a result of H. Liu, \cite[Proposition~5.1]{Liu}, $I_N$ is the small $J$-function of $\mathbb{P}^N$ with twistings given by $\frac{1}{\mathcal{Eu}_{\lambda}(T^*\mathbb{P}^N)}$ and $\operatorname{det}(T^*\mathbb{P}^N)$. Equivalently,
it is the localized stable-map theory of $T^*\mathbb P^N$ with the level -1 structure with respect to the bundle defined by pulling back $T^*\mathbb{P}^N$ from the base. 

\begin{prop}\label{cotangent-difference-equation}
The series $\widehat I_N$ satisfies
\begin{equation}\label{cotangent-q-difference}
 \prod_{j=0}^{N}
 (1-P\Lambda_j^{-1}q^{Q\partial_Q})\widehat I_N
 =Q\prod_{j=0}^{N}
 (1-\lambda P\Lambda_j^{-1}q^{Q\partial_Q})\widehat I_N.
\end{equation}
Its degree-$d$ coefficient after division by $1-q$ has $q$-degree
$-(N+1)d$.
\end{prop}

\begin{proof}
The ratio of consecutive coefficients is
$\prod_{j=0}^{N}(1-yP\Lambda_j^{-1}q^{d-1})/
(1-P\Lambda_j^{-1}q^d)$,
which proves the recurrence.  The numerator contributes
$(N+1)\binom d2$ and the denominator contributes
$(N+1)\binom{d+1}2$, which proves the claim about degrees.
\end{proof}

Quantum triviality and formal Nakayama now give a complete presentation.

\begin{thm}\label{cotangent-ring-comparison}
The stable-map twisted ring defined by $\widehat I_N$ is
\begin{equation*}
 QK_T^{tw}(\mathbb P^N)\cong
 \frac{K_T(\mathrm{pt})[\lambda^{\pm1},P^{\pm1}][[Q]]}
 {\left(
 \prod_{j=0}^{N}(1-P/\Lambda_j)
 -Q\prod_{j=0}^{N}(1-yP/\Lambda_j)
 \right)}.
\end{equation*}
Under $\lambda=\hbar^{-1}$, $Q=(-1)^{N+1}z\hbar^{(N+1)/2}$, and
$P\mapsto s$, its defining relation is equivalent to the Bethe equation
above.  Hence it is isomorphic to the PSZ quasimap quantum
$K$-ring of $T^*\mathbb P^N$, with the $\mathcal{O}(1)\mapsto s$. 
\end{thm}

\begin{proof}
The degree estimate in Proposition~\ref{cotangent-difference-equation}
implies
$\widehat P=P$.
Equation \eqref{cotangent-q-difference} yields the resulting relation. Theorem~\ref{K-Nakayama} shows that it is the full stable-map ring.

For $y=\hbar^{-1}$,
$\prod_{j=0}^{N}(s-\Lambda_j)/(\hbar\Lambda_j-s)
=(-1)^{N+1}\hbar^{-(N+1)}
\prod_{j=0}^{N}(1-s/\Lambda_j)/(1-s/(\hbar\Lambda_j))$.
The stated substitution therefore turns one
relation into the other.  The PSZ presentation uses the quantum tautological
class as its generator, which gives the stated isomorphism.
\end{proof}

\section{Grassmannians}\label{grassmannians}

Let $X=\operatorname{Gr}(k,n)$ be the Grassmannian of $k$-planes in
$\mathbb C^n$, with tautological bundle $\mathcal S$ and quotient bundle
$\mathcal Q=\mathbb C^n/\mathcal S$.  The diagonal torus acts on
$\mathbb C^n$ with characters $\Lambda_1,\dots,\Lambda_n$.

The GIT presentation we use is
$\operatorname{Gr}(k,n)=\operatorname{Hom}(\mathbb C^k,\mathbb C^n)
\sslash GL_k$, with abelian quotient $Y=(\mathbb P^{n-1})^k$.  If $P_i$ is the Hopf line
bundle on the $i$th factor, then
\begin{equation*}
 K_T(Y)=
 \frac{K_T(\mathrm{pt})[P_1^{\pm1},\dots,P_k^{\pm1}]}
 {\left(\prod_{a=1}^n(1-P_i/\Lambda_a):1\leq i\leq k\right)}.
\end{equation*}
The Weyl group $S_k$ permutes the $P_i$, and the classical map sends
$e_r(P_1,\dots,P_k)$ to $\wedge^r\mathcal S$.  The root bundle consists
of the line bundles $P_i/P_j$, $i\neq j$.

Specializing to the symmetrized theory, Givental and Yan prove that the map $\Phi_Q$ gives a correspondence $\mathcal{L}^{tw}_Y$ and $\mathcal{L}_{Gr(k,n}$. 
Givental and Yan compute the equivariant small $J$-function of the
Grassmannian \cite[Theorem~2]{Xiaohan}:
\begin{align*}
 J_X={}&(1-q)\sum_{d_1,\dots,d_k\geq0}Q^{d_1+\cdots+d_k}
 \frac{1}{\displaystyle
 \prod_{i=1}^k\prod_{a=1}^n\prod_{m=1}^{d_i}
 (1-q^mP_i/\Lambda_a)}\notag\\
 &\quad\times
 \prod_{i,j=1}^k
 \frac{\prod_{m=-\infty}^{d_i-d_j}(1-q^mP_i/P_j)}
 {\prod_{m=-\infty}^{0}(1-q^mP_i/P_j)}.
\end{align*}
Products with lower limit $-\infty$ are ratios, and the $i=j$ factors are
$1$.

It is the image under $\Phi_Q$ of the following value of the twisted $J$-function of $Y$:
\begin{align*}
 I_Y^{\lambda}:={}&(1-q)
 \sum_{d_1,\dots,d_k\geq0}\prod_{i=1}^kQ_i^{d_i}
 \frac{1}{\displaystyle
 \prod_{i=1}^k\prod_{a=1}^n\prod_{m=1}^{d_i}
 (1-q^mP_i/\Lambda_a)}\notag\\
 &\quad\times
 \prod_{i\neq j}
 \frac{\prod_{m=-\infty}^{d_i-d_j}
 (1-\lambda q^mP_i/P_j)}
 {\prod_{m=-\infty}^{0}(1-\lambda q^mP_i/P_j)}.
\end{align*}

Since the $d=0$ term of this is $1-q$, the input $t$ lies in $\Lambda_+$.

\subsection{Degree bounds}

Write $I_Y^\lambda/(1-q)
=\sum_{\mathbf d\geq0}Q_1^{d_1}\cdots Q_k^{d_k}I_{\mathbf d}(q)$.
Rearranging the root factors as in \cite[Remark~2]{Xiaohan} shows that the excess of
denominator degree over numerator degree is
\begin{equation*}
 D(\mathbf d)
 :=n\sum_{i=1}^k\binom{d_i+1}{2}
 -\sum_{d_i>d_j}\binom{d_i-d_j+1}{2},
 \deg_q I_{\mathbf d}=-D(\mathbf d).
\end{equation*}

\begin{lem}\label{grass-degree-bound}
For every nonzero $\mathbf d\in\mathbb Z_{\geq0}^k$,
$D(\mathbf d)\geq(n-k+1)\sum_{i=1}^kd_i$.  Consequently
$\deg_q I_{\mathbf d}< -\sum_i d_i$ and
$\deg_q I_{\mathbf d}<-(n-k)\max_i d_i$.
\end{lem}

\begin{proof}
Put $B_i=\binom{d_i+1}{2}$.  Whenever $d_i>d_j$,
$\binom{d_i-d_j+1}{2}\leq B_i-B_j\leq B_i$.
For each $i$ there are at most $k-1$ smaller entries, and therefore the
second sum in the definition of $D(\mathbf d)$ is at most
$(k-1)\sum_iB_i$.  Hence
$D(\mathbf d)\geq(n-k+1)\sum_iB_i
\geq(n-k+1)\sum_i d_i$.
Since $n-k+1\geq2$ and $\mathbf d\neq0$, both displayed strict
inequalities follow.
\end{proof}

The above degree bounds also imply the input $t(q)$ vanishes at $q=\infty$. Since it is $W$-invariant and $\Phi_Q(t(q))=0$, because $I$ specializes to the small $J$-function, we can
apply Corollaries~\ref{coras}--\ref{corbs}, which yield:
\begin{align}
 \Phi_Q(e_r^{\star_Y}(P))\star_Z
 &=\bigl(\wedge^r\mathcal S\bigr)\star_Z,
 &&0\leq r\leq k,
 \label{grass-elementary-triviality}\\
 \Phi_Q(h_r^{\star_Y}(P))\star_Z
 &=h_r(\mathcal S)\star_Z,
 &&0\leq r\leq n-k,
 \label{grass-complete-triviality}
\end{align}
where $f^{\star_Y}(P)=f(P_1\star_Y,\dots,P_k\star_Y)1$ is the quantum
evaluation from Corollary~\ref{corbs}, and $h_r$ denotes the complete
homogeneous symmetric polynomial.

\subsection{Difference equations and the Bethe Ansatz}

For a finite alphabet $A=(A_1,\dots,A_m)$, write
$e_r(A):=\sum_{1\leq i_1<\cdots<i_r\leq m}A_{i_1}\cdots A_{i_r}$,
with $e_0(A)=1$ and $e_r(A)=0$ for $r<0$ or $r>m$.  Define
$E_A(y):=\sum_{r=0}^m e_r(A)y^r=\prod_{i=1}^m(1+yA_i)$.
We also use the complete homogeneous functions determined by
$\sum_{r\geq0}h_r(A)y^r=E_A(-y)^{-1}
=\prod_{i=1}^m(1-yA_i)^{-1}$.
Thus $P$ denotes the alphabet $(P_1,\dots,P_k)$, while the symbol $\Lambda$
inside $e_r(\Lambda)$ or $E_\Lambda(y)$ denotes the alphabet of equivariant
characters $(\Lambda_1,\dots,\Lambda_n)$, not the coefficient algebra
$\Lambda$ used elsewhere in the paper.

For $D_i=Q_i\partial_{Q_i}$, direct comparison of adjacent coefficients in
the definition of $I_Y^\lambda$ gives
\begin{align*}
 &\prod_{j\neq i}
 \left(1-\lambda q\frac{P_j}{P_i}q^{D_j-D_i}\right)
 \prod_{a=1}^n\left(1-\frac{P_i}{\Lambda_a}q^{D_i}\right)I_Y^\lambda
 \notag\\
 &=Q_i\prod_{j\neq i}
 \left(1-\lambda q\frac{P_i}{P_j}q^{D_i-D_j}\right)I_Y^\lambda,
 \quad i=1,\dots,k.
\end{align*}
Applying Theorem~\ref{abelianized-relations}, then setting
$q=1$, $Q_i=Q$, and $\lambda=1$, yields the following relations in $QK^{tw}(Y)$, up to the ideal $\mathscr{I}:$.  
\begin{align}
 \prod_{a=1}^n(1-P_i/\Lambda_a)
 &=Q\prod_{j\neq i}
 \frac{1-P_i/P_j}{1-P_j/P_i}\notag\\
 &=(-1)^{k-1}Q\frac{P_i^k}{e_k(P)},
 \quad i=1,\dots,k.
 \label{grass-Bethe-relations}
\end{align}
 We remark that quations \eqref{grass-Bethe-relations} are the
Bethe-Ansatz equations for
$QK_T(\operatorname{Gr}(k,n))$
appearing in the Coulomb branch of the firled theory described in \cite[Theorem~1.2]{gu2022quantum}.

Symmetric combinations of equations 8.3 descend as relations via $\Phi_Q$. We describe below an approach to obtain convenient such combinations:

Define the degree-$n$ polynomial
\begin{equation*}
 F(\xi):=\prod_{a=1}^n(\Lambda_a-\xi)
 -(-1)^{k-1}Q\frac{e_n(\Lambda)}{e_k(P)}\xi^k.
\end{equation*}
Then \eqref{grass-Bethe-relations} says $F(P_i)=0$.  In a splitting
extension, write the full multiset of roots of $F$ as
$w=P\sqcup\overline P$, where $P=\{P_1,\dots,P_k\}$ and
$|\overline P|=n-k$.
Vieta's formulas give
\begin{equation*}
 e_r(w)=
 \begin{cases}
 e_r(\Lambda),&r\neq n-k,\\[2mm]
 e_{n-k}(\Lambda)+Q\dfrac{e_n(\Lambda)}{e_k(P)},&r=n-k,
 \end{cases}
\end{equation*}
Equivalently,
\begin{equation}
\label{generatingv}
 E_P(y)E_{\overline P}(y)
 =E_\Lambda(y)+Qe_n(\Lambda)e_k(P)^{-1}y^{n-k}.
\end{equation}
Since the top Vieta relation also gives
$e_n(\Lambda)/e_k(P)=e_{n-k}(\overline P)$, this becomes
\begin{equation}
\label{generatingv2}
E_P(y)E_{\overline P}(y)
=E_\Lambda(y)+Qe_{n-k}(\overline P)y^{n-k}.
\end{equation}
Expanding $E_\Lambda/E_P$ expresses the coefficients in terms of
$e_a(\Lambda)$ and $h_b(P)$.  The identities
\eqref{grass-complete-triviality} therefore imply
\begin{equation*}
 \Phi_Q(e_r(\overline P))=
 \begin{cases}
 \wedge^r\mathcal Q,&0\leq r<n-k,\\[2mm]
 \dfrac{1}{1-Q}\det\mathcal Q,&r=n-k.
 \end{cases}
\end{equation*}
Indeed, the coefficient of $y^{n-k}$ in
\eqref{generatingv2}, after division by $E_P(y)$, reads
$(1-Q)e_{n-k}(\overline P)
=[y^{n-k}]E_\Lambda(y)/E_P(y)$,
and the right-hand side specializes to $\det\mathcal Q$.

The identity obtained below is the equivariant quantum Whitney relation of
Gu--Mihalcea--Sharpe--Zou \cite[Theorem~1.1]{gu2022quantum}.

\begin{thm}[Quantum Whitney relation]\label{grass-quantum-Whitney}
In the small equivariant quantum $K$-ring of the Grassmannian,
\begin{equation}\label{eq:gr:lambda}
 \lambda_y(\mathcal S)\star\lambda_y(\mathcal Q)
 =\lambda_y(\mathbb C^n)
 -\frac{Q}{1-Q}y^{n-k}
 \bigl(\lambda_y(\mathcal S)-1\bigr)\star\det\mathcal Q.
\end{equation}
Equivalently, with $\wedge^m\mathcal S=0$ for $m\notin[0,k]$,
\begin{align}\label{grass-Whitney-components}
 \sum_{a+b=\ell}\wedge^a\mathcal S\star\wedge^b\mathcal Q
 ={}&\wedge^\ell\mathbb C^n\notag\\
 &-\frac{Q}{1-Q}\det\mathcal Q\star
 \left(\wedge^{\ell-n+k}\mathcal S-\delta_{\ell,n-k}\right).
\end{align}
\end{thm}

\begin{proof}
Apply $\Phi_Q$ to \eqref{generatingv2}.  By
\eqref{grass-elementary-triviality} and the displayed specialization of
$e_r(\overline P)$,
\begin{equation*}
 \Phi_Q(E_P)=\lambda_y(\mathcal S),
 \Phi_Q(E_{\overline P})
 =\lambda_y(\mathcal Q)
 +\frac{Q}{1-Q}y^{n-k}\det\mathcal Q.
\end{equation*}
\end{proof}

These relations yield the following presentation of $QK_T(Gr(k,n))$:

\begin{cor}[Whitney presentation]\label{grass-Whitney-presentation}
The coefficients of \eqref{eq:gr:lambda}, with generators
$\wedge^i\mathcal S$ for $1\leq i\leq k$ and
$\wedge^j\mathcal Q$ for $1\leq j\leq n-k$, give a presentation of
$QK_T(\operatorname{Gr}(k,n))$ over $K_T(\mathrm{pt})[[Q]]$.
\end{cor}

\begin{proof}
At $Q=0$ the relation is the classical Whitney identity and its quotient has
the Schur basis indexed by partitions in the $k\times(n-k)$ rectangle.  The
quantum relations preserve this finite free basis (equivalently, one uses
the straightening argument in the proof of
\cite[Theorem~1.2]{gu2022quantum}).  The map from the
presented algebra to quantum $K$-theory is surjective and is an isomorphism
modulo $Q$, so Theorem~\ref{K-Nakayama} proves that it is an isomorphism.
\end{proof}

\section*{Acknowledgements}
The author is thankful to Eric Sharpe, Leonardo Mihalcea, Weihong Xu, Xiaohan Yan, Kamyar Amini, and Qasim Shafi for their discussions regarding the content of this work. The author is thankful to Leo Herr for feedback on readability.
\clearpage
\printbibliography

\end{document}